\documentclass[12pt]{article}
\usepackage{graphicx} 
\usepackage{amscd}
\usepackage{amsmath}
\usepackage{amssymb}
\usepackage{amsthm}
\usepackage{amstext}
\usepackage{amsopn}
\usepackage{texdraw}
\usepackage{graphics}
\usepackage{mathptmx}       
\usepackage{helvet}         
\usepackage{courier}        
\usepackage{makeidx}         
\usepackage{graphicx}        
\usepackage{multicol}        
\usepackage[bottom]{footmisc}
\usepackage{braket}
\usepackage{multirow,array}
\newtheorem{theorem}{Theorem}[section]

\theoremstyle{definition}
\newtheorem{definition}[theorem]{Definition}

\theoremstyle{remark}
\newtheorem{remark}[theorem]{Remark}
\numberwithin{equation}{section}

\usepackage{mathptmx}       
\usepackage{helvet}         
\usepackage{courier}        
\usepackage{makeidx}         
\usepackage{graphicx}        
\usepackage{multicol}        
\usepackage{multirow}
\usepackage[bottom]{footmisc}
\usepackage{tikz}
\makeindex             

\usepackage{amscd,amsmath,amssymb}

\title{Compendium of Advances in Game Theory: Classical, Differential, Algorithmic, Non-Archimedean and Quantum Games}
\author{Bourama Toni\\
Department of Mathematics\\
Howard University, Washington DC 20059\\
\textcolor{blue}{E-mail: bourama.toni@howard.edu}}

\date{April 2025}

\begin{document}



\maketitle

\begin{abstract}
This compendium features advances in Game Theory, to include:
\begin{itemize}
\item Classical Game Theory: Cooperative and non-cooperative. Zero-sum and non-zero sum games. Potential and Congestion games. Mean Field games. Nash Equilibrium, Correlated Nash Equilibrium and Approximate Nash Equilibrium. Evolutionary Game Theory
\item Intelligent Game: Differential Game Theory. Algorithm Game Theory and Security Games.
\item Quantum Games and Quantumization of classical games such as the Battle of the Sexes
\item Non-Archimedean and p-adic game theory and its growing relevancy as the domains of game-theoretic application expands. 
\item p-adic quantum game to leverage and combine the distinguishing features of non-Archimedean analysis and quantum information theory. This is a novel game-theoretic approach with great potential of application.
\end{itemize}
In times of exponential growth of artificial intelligence and machine learning and the dawn of post-human mathematical creativity, this compendium is meant to be a reference of choice for all game theory researchers.
\end{abstract}

 Complex Hilbert space. p-adic Hilbert space. p-adic quantum game.

\newpage

 {\bf In Memoriam: John F. Nash (1928-2015)}
 
\begin{figure}[h!]
\centering
\includegraphics[height=9cm,width=.90\textwidth]{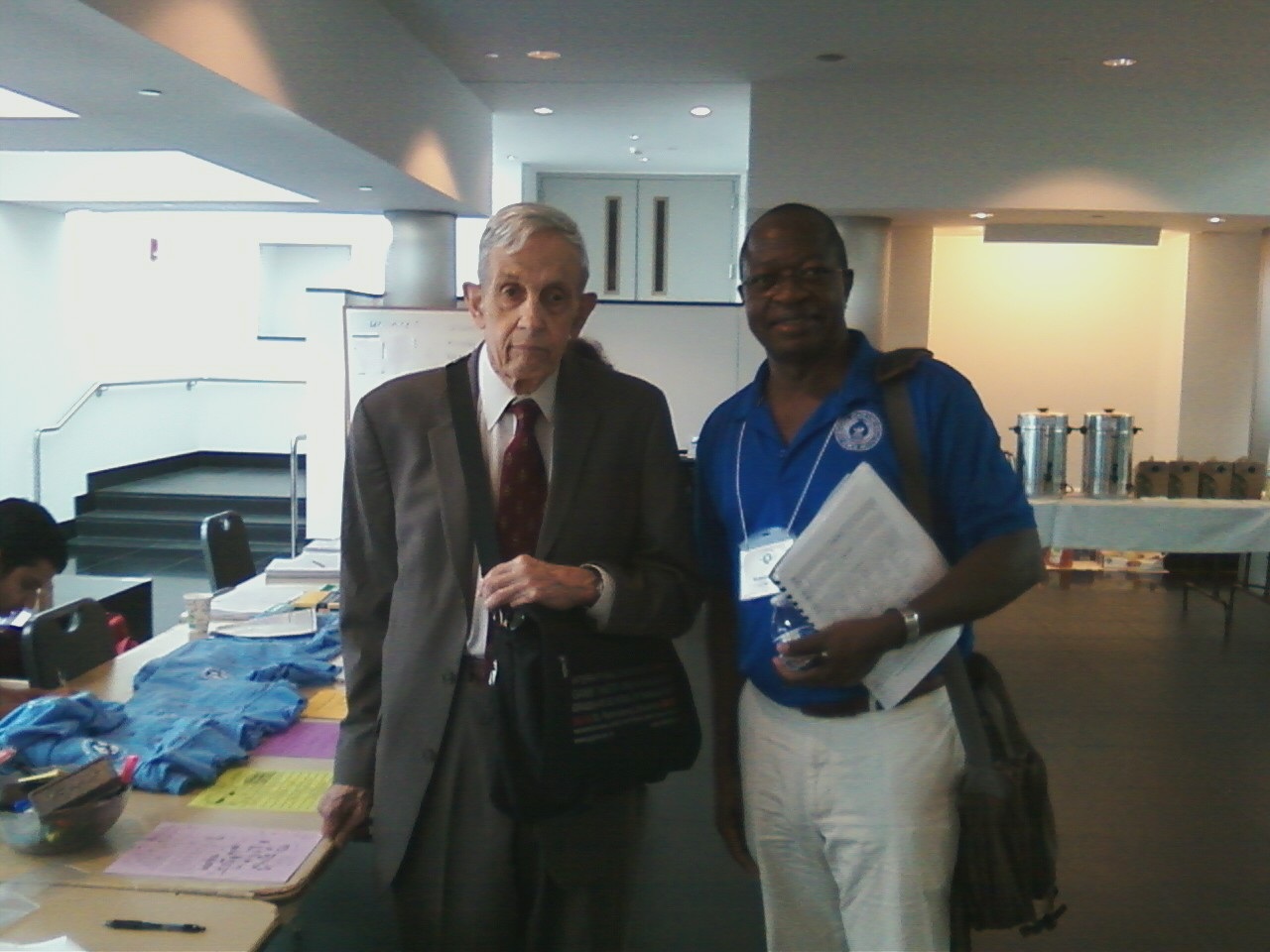}

{\bf Towards a concept of Nash Limit Cycle, Stony Book July 2013}
\end{figure}

At the July 2013 International Conference on Game Theory at Stony Brook,  we were fortunate to meet and discuss with the late John Nash some of our game-theoretic ideas, in particular, the concept of \textcolor{blue}{Nash Limit cycles} and its application to socio-cultural evolution. His inputs were valuable to us then, as they are today, and we dedicate this work and what it is worth to his memory and his everlasting contributions to this old age Theory of Games. It has been a privilege to get to know John Nash. His genuine kindness and keen show of interest in our game-theoretical work and the STEAM-H series is fondly remembered: he so  kindly autographed our first two early books.

\newpage


{\bf OUTLINE}

\begin{enumerate}
    \item Introduction and Background
    \item Evolutionary Game Theory
      \item Intelligent Game
    \item Differential Game Theory
    \item Algorithmic Game Theory and Security Games
    \item Quantum Game Theory
    \item Non-Archimedean Game Theory
    \item p-Adic Quantum Game Theory
    \item A socio-cultural game model
    \item Concluding Remarks
    \item References
\end{enumerate}

\maketitle

\section{Introduction and Background}

\subsection{Raison d'\^etre}

Game Theory has evolved from an intuitive to a formal analysis and understanding of interacting entities, be it, organism, animals, programs, social norms, etc. Its underlying ideas and {\it raison d'\^etre} transpired throughout history: from the biblical game-theoretic story of the  Babylonian {\it Talmud} \footnote{The Talmud, more precisely the Babylonian Talmud, comprises of The Mishna (c. 200 CE), a written compendium of Judaism's Oral Law and The Gemara (c. 500 CE), which is a record of discussions by rabbis about the Mishna. It is subdivided into 60 books with the first printed version appearing around 1520. A modern English translation has 73 volumes. The Mishna is in Hebrew, but the Gemara is Aramaic; The Talmud is thought to be largely incomprehensible without the commentary of the French rabbi Rashi.} {\it Estate Division Problem} 
\footnote{According to the Babylonian Talmud, a man dies without paying all his three creditors but leaves an estate too small to pay his debts. Creditor 1 is owed 100, Creditor 2 is owed 200 and Creditor 3 is owed 300. If the estate is 100, each creditor gets 33 1/3, If the estate is 200, Creditor 1 gets 50, Creditors 2 and 3 get 75 each, And finally if the estate is 300, Creditor 1 gets 50, Creditor 2 gets 100, and Creditor 3 gets 150. The estate division algorithm is based on the legal system's principles and precedents. The problem was completely solved by the mathematician Robert Aumann and Michael Maschler in the 1980s at the Hebrew University of Jerusalem.(Brams 1980 and Aumann and Mashler, 1985, Schecter 2012).}  to the works of Descartes, \footnote{ Descartes' {\it Cogito ergo sum} or "I Think therefore I am" is connected to the idea of rational decision making central to Game Theory: an individual's choice based on own understanding and assessment of a situation.}
Darwin\footnote{Evolution \'a la Darwin has been combined with Game Theory to explore and understand how organisms evolve and adapt.}, 
and the teachings of the Chinese warrior philosopher Sun Tzu (Sun Tzu, 1988). \footnote{Sun Tzu, circa 544-496 BCE is mostly known for The Art of War stating:

{\it Knowing the other and knowing oneself, In one hundred battle no danger

Not knowing the other and knowing oneself, One victory for one loss

Not knowing the other and not knowing oneself, In every battle certain defeat}

(S. Tzu in the Denma Translation, Shambala Library 2002)}. 

Most societies, past and present, have been and are fascinated by game playing and gambling, in particular, {\it games of chance or of strategies}. \footnote{Playing cards appeared in Europe around the tenth century. A precursor of the six-faced die, called {\it Astragalus},  was found on Sumerian, Assyrian and Egyptian archaeological sites, circa 3600 BCE. Lotteries are from the first century Roman Empire. } (i.e., outcomes determined by chance or by skills and abilities): Dice, Cards, Roulette, Poker, Craps, Lotteries, Spin-to-win, Random Prize draws. Useful and practical examples for experiments, games and gambling have also contributed to the development of the mathematical theory of Probability. Indeed, there has been some serious scientific excitements and achievements strongly related to Games and Gambling; notably, the work of Girolamo Cardano in his {\it The Book on Games of Chance} around 1520, who also excelled in Chess and Backgammon and other dice games; Blaise Pascal and Pierre de Fermat have made foundational contributions to probability and gambling questions around 1654, encouraged by their compulsive gentleman gambler friend, Antoine Gombaud, a.k.a the Chevalier de M\'er\'e\footnote{Self-knighted Chevalier de M\'er\'e, he is a typical french {\it salon theorist} known for his essays "L'honn\^ete homme'' and "Discours de la vraie honn\^etet\'e" }. This has later led to the famous {\it Pascal's wager}, considered a game-theoretic approach to the belief in God: that is, it is "probabilistically prudent" to believe in God. (See Packel 2006). More notably,   Christiaan Huygens around 1657 introduced the concept of expectation in his treatise {\it On reasoning in Games of Chance (De Ratiociniis in ludo aleae},

What exactly is game theory? In essence, Game Theory\footnote{Game Theory has been popularized by various means and is the theme of the 2001 Best Picture and Best Director, {\it Beautiful Mind}, featuring Russell Crowe playing John Nash, the nerdy Princeton student who's, alledgedly, life and dating experience led him to the Nash Equilibrium solution concept that ultimately ended in a 1994 Noble Price Award in Economics. (Nasar, 1998)} is a mathematical framework (American Mathematical Classification code 91A) that allows analysis of strategic decision making between multiple entities called {\it the Players}. It has all the features of a mathematical construct: a concise set of concepts and assumptions, many fundamental theorems, and many real-life applications where human and social behaviors are mathematically analyzed. It has traditionally carries a {\it rationality assumption}, however, more concerned with the "pursuit of happiness" by the decision-makers (i.e., not necessarily altruists).

Game theory or Interactive decision-making was first initiated in the field of microeconomics. See {\it Theory of Games and Economic Behavior} by the mathematician John Von Neumann\footnote{ Von Neumann's 1928 paper founded the theory of two-person zero-sum games, making him the mathematician most closely related to the creation of the theory of games. However, Emile Borel preceded him in formulating a theory of games.} and the economist Oscar Morgenstern. The book has a seismic effect on quantitative social sciences (economics, political science and psychology), providing the long-awaited mathematical equivalency of mathematical physics.

Interactive strategic decision-making schemes prevail across most disciplines, leading to the expansion of the applicability of game-theoretic methodology over the years to areas such as Stock Markets and Financial Management, Military War Games, National Security and Anti-terrorism, Bridge games, Hunting Party, Politics, Evolutionary Biology, Computer Science. The {\it actions/strategies} of the {\it players} should result in the {\it best possible} consequences or outcomes according to their preferences. (Axelrod, 1984. Bisin et al., 2001. Bowles et al. 2004).

First, let us agree on what a decision-making is. A good related reading is found in the book "Game Theory" by James Webb. (Webb, 2007)

\subsubsection{ Optimum and Rational Decision Making}

In a game play, {\it Decision making} actually refers to making \textcolor{blue}{optimum decision-making}. That is, from an available set of courses of action or possible behaviors one first determines the outcome of each one of these actions, and then chooses the {\it preferred} one, in the sense of the one that maximizes a personal utility or payoff, numerical or otherwise (e.g. something of great value to one, be it monetary, social norms, etc.). Maximizing/minimizing easily translates in a mathematical framework (extremum of a function in basic calculus). Therefore, decision-making sets in motion a course of action that leads to an outcome, which is certain in the absence of randomness.  A \textcolor{blue}{payoff function} is then determined to associate a numerical value with each action. The decision-maker strives to maximize/optimize such a payoff. Note that a payoff is "invariant" under certain changes, the so-called {\it affine transformation or rule} assigning the same value of a payoff $\pi(x)$ to its affine transformed payoff $a\pi(x)+b.$ In general, payoffs model individual preferences during the social interactions. The so-called {\it Darwinian fitness} in Evolutionary game models. 

\textcolor{blue}{Uncertainty} is an inherent part of any decision-making. We then talk about \textcolor{blue}{expected outcome/expected payoff}, to be represented as a random variable  and its associated probability distribution. For example, decision-making with respect to the performance of the stock market. An uncertain outcome is also called {\it a lottery,} which implies a finite probability distribution over the set of all payoffs.  Of course, randomization is not required to maximize one's payoff; but it can serve as a tie-breaker as seen below in the case of multiple acceptable outcomes or equilibria. We will refer to a plan of action as a {\it strategy}, {\it pure} in the absence of randomization, and {\it mixed} otherwise. 

\subsection{Variation on the nature of Games}

A player could be any decision-making entity such as an automaton, a machine, a program, a person or an animal, a living cell or a molecule. The term {\it Game} refers to in any interactive situation leading to the sharing of clearly quantifiable benefits, the {\it payoffs/utilities} and costs. Traditionally, a fundamental premise has been the assumed {\it rationality} of all players. The strategies may be tightly coupled, strongly correlated (i.e., {\it entangled} as in quantum strategies), allow {\it superposition} (e.g., quantum game) or may allow a probability distribution for a random selection ({\it mixed strategies}). Strategies termed {\it pure} refer to simple and consistent action plans specific to a given player (acting with certainty) during a game (e.g., playing {\it rock} in the popular {\it rock-paper-scissors} game). The continuity of probabilities ensure an infinitely many mixed strategies for the player. A {\it totally mixed strategy} implies the assignment of strictly positive probability to every pure strategy (also called a {\it degenerate} mixed strategy with probability one).

Along with the rationality assumption for every player, goes also the {\it selfishness} assumption. That is, players strive to maximize their self-interest, i.e., the individual payoffs during a game. Preferences in game theory may include altruistic motivations, moral principles, social constraints. Oftentimes, any perceived cooperative behavior is transient and driven by selfish objectives. However, we distinguish such game called a \textcolor{blue}{noncooperative/competitive game from cooperative or coalitional game}. In the latter game, players form coalitions or groups, with sometimes external enforcement, and focusing on, e.g., surplus or profit-sharing among the coalitions.

Competitive games may also describe the so-called \textcolor{blue}{zero-sum games} in which the sum of the payoffs of all players is zero regardless of their actions/moves/strategies: in its two-person version one player's loss corresponds to another player's gain. A so-called payoff matrix is often used to represent a two-player zero-sum game. In solving such a game, one seeks a set of strategies to minimize the maximum loss, the {\it Minimax Theorem} or, inversely to maximize the minimum payoff, the {\it Maximin Theorem}. The solution or outcome is called a {\it pure saddle point}. The antithesis of such games is the {\it common interest game}: players have perfectly aligned interests.
A particular class of zero-sum games or win-loss game scenarios is that of \textcolor{blue}{Combinatorial Games.}\footnote{Initiated around 1901 by Charles Leonard Bouton, the theory of combinatorial game was further developed by Roland Percival Sprague and Patrick Michael Grundy in the 1930s, continued with John Milnor and Olof Hanner in the 1950s.} in which players alternate in moving. This class of games is also benefiting from the blooming development of artificial intelligence and deep neural networks. It is typically a two-player game with no hidden information (i.e., full disclosure of the game's position) and no chance elements (i.e., the next position of the game is totally determined by a player's move). There are in general two winning conditions along the "divide and conquer" approach: "Last move wins" called the {\it Normal play} or "last move loses" called the {\it Mis\`ere play}. See {\it Winning Ways for your Mathematical Plays}, Berlekamp et al., 1982. Combinatorial Game is rooted in ancient board games (e..g, Go, Chess, Checkers) and in the so-called {\it recreational mathematics} such as {\it The Tower of Hanoi}.\footnote{Francois-Edouard Anatole Lucas, 1842-1891, is seen as the greatest French recreational mathematician and invented the Tower of Hanoi game in 1883 and whose original is the Conservatoire National des Arts et M\'etiers-Mus\'ee National des Techniques. $n$ discs require $2^n-1$ moves to solve the problem}. Indeed, mathematically the Tower of Hanoi corresponds to a {\it Hamiltonian Circuit} on the n-dimensional cube. Other well-known recreational mathematics include the {\it Chinese Rings}, aka {\it Baguenaudier, Cardan's suspension, Cardano's rings, Devil's needle, five pillars puzzle} and the Chinese {\it I-Ching Hexagrams}\footnote{Book of Changes with 64 hexagrams corresponding to nature fundamental forces: the Yin for the 6 broken lines and the Yang for the 6 unbroken ones.}. Well-known Combinatorial Games include Chess, Checkers, Tic-Tac-Toe, Go and Connect Four. Well-known games but not considered combinatorial are Bridge, Backgammon, Poker, Snakes and Ladders.

Some competitive games such as Chess may feature {\it perfect information} to refer to the fact that each player is perfectly aware of the previous actions of all other players, including of course the initial state of the game in a sequential game (with predefined order). On the other hand the game may be an {\it imperfect information} one, with some players not accessing the entirety of other players' actions.

In classical game theory, there are typically three dominant mathematical forms(1) Normal Form or Strategic Form Game (SFG) (2) Extensive Form Game (EFG) (3) Coalition Form Game (CFG) dealing with options for subsets of players. The SFG consists of players making the relevant decisions with the strategies available, the payoffs are the rewards contingent upon the actions of all the players in the game. The EFG places the emphasis on the timing of the decisions to be made, as well as the information available, and is representable by a {\it decision or game tree}. See more details in Fudenberg et al. 1991; Webb 2007.

\subsubsection{Normal or Standard Form Game}

 The most common form of classical games is the so-called \textcolor{blue}{Normal Form/Standard Form}:  all players make decision simultaneously without knowledge of other players' actions. A player controls only their own actions, and is interested only in action profiles that maximize their own payoff/utility. 

A formal definition or representation of a {\it Normal Form} game is as follows, using the terms {\it action} and {\it strategy} interchangeably, as well as {\it utility} and {\it payoff.}

\begin{definition}
A {\it Strategic Normal Form Game} is described by a triplet $\mathcal G=\langle \mathbf P,(\mathbf S_i)_{i\in\mathbf I_n},(u_i)_{i\in\mathbf I_n}\rangle ,$ where $\mathbf I_n=\{1,\ldots,n\},$ and 
\begin{enumerate}
  \item $\mathbf P$ is a finite set of the $n$ players $P_{i=1,\dots,n}$
  \item $\mathbf S_i$ is the non-empty set of {\it pure strategies} or available actions for player $P_i.$ The set $S_i$ can be written as $S_i=\{s_1^1, s_i^2,\ldots,s_i^{\alpha},\ldots,s_i^{n_i}\},$ where the $s_i^{\alpha}$ are the {\it pure strategies} player $p_i$ has access to, and $n_i$ is the finite number of these pure strategies. (Traditionally the payoff functions are real-valued.)
  
  \item $u_i:\mathbb S=\prod \mathbf S_{i\in \mathbf I_n}\rightarrow \mathbb F$ is the {\it utility} or {\it payoff} function for player $P_i.$ Thus a tuple of payoff functions $u=(u_1, \ldots,u_n).$ For $u=(u_1,\ldots,u_n)$ the set of payoffs is given by $u(S)=(u_i(S))_{i=1,\ldots,n}.$
  
  \item The vector $s=(s_1,\ldots,s_n)\in \mathbb S$ denotes the {\it pure strategy/ action profile} or {\it outcome}
  \item We denote by $s_{-i}=(s_1,\ldots,s_{i-1},s_{i+1},\ldots,s_n)$ the vector of strategies taken by all other players except player $P_i.$ That is, $s=(s_i,s_{-i}).$ And accordingly $\mathbb{S}_{-i}=\prod_{j\ne i}\mathbf S_j$ describes all action profiles for all players except for $P_i.$
  \end{enumerate}
  \end{definition}

 A "play" of the game refers to the n-tuple strategy profile or pure state 
 $$
 G:=s^{\alpha}=(s_1^{\alpha_1},\ldots,s_n^{\alpha_n})
 $$
where player $P_i$ chooses to play the strategy $s_i^{\alpha_i} \in S_i,$ and received a reward or payoff computed by the payoff function 
$$ 
u_i(s^{\alpha})= u_i(s_1^{\alpha_1},\ldots,s_n^{\alpha_n})
$$
The reward or payoff, in general, depends on the strategies played by all the players, under the rational assumption that each tries to maximize their payoff. Thus a need of a compromise or an {\it equilibrium} or {\it solution} of the game.

For games of "complete information" players know the strategy and payoff spaces of each other.

\begin{remark}
Note that, classically, the utility  or payoff function is real-valued. The co-domain $\mathbb R$ is endowed with the Euclidean/Archimedean norm/absolute value $|.|_{\infty}$ . A recent advance, described below, is to consider p-adic/non-Archimedean valued payoff function to more effectively address nuances and the possible hierarchy in the payoffs structures. That is, the co-domain is the payoff function $u_i$ is $\mathbb{Q}_p,$ the set of p-adic numbers defined below, endowed with p-adic norm $|.|_p.$ See Gouv\^ea 1997.
\end{remark}

\begin{definition} (Dominant Strategy)
A pure strategy/action or value $\alpha\in S_i$ is said to be {\it (strictly) dominated} by a strategy $\beta\in S_i$, $\forall i=1,\ldots,n$ denoted by $ \alpha <\beta$ when $u_i(\alpha,s_{-i})<u_i(\beta,s_{-i}).$
\end{definition}

\begin{remark}
Dominated strategies/actions/profiles are {\it rationally unjustifiable}. However removing them in a given game may lead to the emergence of other dominated strategies in the subgame, which results inductively in the so-called {\it iteratively dominated strategies}. {\it Iteratively undominated strategies} survive all rounds of elimination by dominance, resulting sometimes in instances of games called {\it dominance solvable} where the undominated strategies form a singleton.(see Fisher 1922, Fudenberg and Tirole 1991).
\end{remark}

\subsubsection{Solution Concepts}

Most of the time games cannot be solved just by removing dominated strategies. The {\it solution concept or equilibrium} refers to a balance of players' strategies in such a way that no player has a motive to change unilaterally. The equilibrium is {\it strong} if it does not allow any group deviation, and it is {\it Pareto efficient or optimal}\footnote{Vilfredo Pareto, Italian mathematician, defines the efficiency as "A situation where there is no way to rearrange things to make at least one person better off without making anyone else worse off".} if there is no other outcome to make all players better off. The rationality assumption indeed forces the appearance of an equilibrium solution.

We then refer to the {fundamental solution } concept that is needed to capture the steady state of play with each player acting optimally given a "conjecture" about the behavior of other players. It is the {\it Nash Equilibrium} of the game to characterize strategy profiles resilient to unilateral deviations. It is akin to a compromise in a game where every player endeavors to maximize their own selfish payoff. Formally

\begin{definition}

A strategy profile $s^*=(s^*_1,\ldots,s^*_n)=(s^*_i,s^*_{-i})\in \mathbb S $ is a \textcolor{blue}{Nash Equilibrium (NE)} if
$$
u_i(s^*)=u_i(s^*_i,s^*_{-i})\ge u_i(s_i,s*_{-i}),\quad \forall s_i\in S_i.
\eqno(3.1)
$$
\end{definition}

\begin{remark}
\noindent
\begin{itemize}

  \item The {\it Nash Equilibrium} must be {\it strategically stable} to unilateral deviations: No player can unilaterally profitably deviate. It is {\it strong} if it is stable to group deviations, as to mitigate mutation and the impact of mutant players.
  
  \item The {\it Nash Equilibrium} exists for a broad class of games as proved by the Nash's famous theorem below. See Nash 1950, 1951, 1953.
  
  \item Some instances of games  may be described by admissible mechanisms for coordination to ensure emergence of desirable collective behavior with respect to a given objective, thus leading to a Nash Equilibrium.  
  
  For example, consider the simple game of  the {\it commuting time to get home from work} which depends not only on the chosen route but also on decisions taken by other drivers. A solution concept or a {\it Nash Equilibrium (NE)} is: {\it In traffic everyone is driving on the right. No single driver (rational) has an interest in driving on the left.}
  \item A NE may result from repeated interactions, even in a game with partial information and primitive decision-making as in the {\it Harper's 33 Ducks Experiment} we here recall. (Harper 1982)

  In the botanic garden of Cambridge University, UK, the biologist David harper experimented with a flock of 33 ducks and 2 bread tossers, 20 meters apart feeding the ducks at regular intervals. Tosser 1 (the least profitable site) has a frequency of supply of 12 pieces of bread per minute  whereas Tosser 2 (the most profitable site) has a frequency of 24 pieces per minute. After about a minute, the number of ducks at the least profitable site stabilizes around 11, with thus 22 ducks at the most profitable site, which indicates a Nash equilibrium has been reached: Switching unilaterally to another site for a duck would result in getting less food. Indeed the amount of food a duck gets does not depend on its choice but also on others' choice. During a transient period, behaving as a typical optimizer most of the ducks rushed to the most profitable site, for some to realize they could get more food at the least profitable site, and therefore switched site. After switching a few times the ducks stick to a given choice. Repeated interactions has led to a Nash Equilibrium in a game where the entities, the duck, could hardly qualified as rational players. Such a process has been termed "{\it t\^atonnement}. (Iterative adjustments in strategies until equilibrium is reached). (Cournot 1838)
  \end{itemize}
\end{remark}

\begin{remark}[Historical Comments]
\noindent
\begin{itemize}
\item An intuitive explanation of the concept of {\it Nash equilibrium} could be traced back to Antoine Augustin Cournot in his {\it Recherches into the Principles of the Theory of Wealth in 1838} in which one could find as well an evolutionary or dynamic idea of the best response correspondence. From Cournot's works Francis Ysidro Edgeworth in his {\it Mathematical Psychics} derived the concept of competitive equilibria in a two-players economy. Then Emile Borel, in {\it Alg\`ebre et calcul des probabilit\'es} published in Comptes Rendus de l'Acad\'emie des Sciences, vol. 185(1927), dealt with mixed strategies, probability distributions leading to stable game.

\item However the modern analysis, a milestone in the history of game theory, is accepted to have been initiated by John von Neumann (1928) and Oskar Morgenstern in their book {\it Theory of Games and Economic Behavior}, drawing from which John Nash provides us with the modern methodological framework. Von Neumann proved there his famous \textcolor{blue}{Minimax Theorem} for zero-sum games.

Other seminal works include that of Nash on \textcolor{blue}{Nash Equilibrium} (Nash, 1951) and on \textcolor{blue}{Nash Barganining} (Nash, 1950), and Shapley on the \textcolor{blue}{Shapley Value} and games with transferable utility (Shapley, 1953, 1967). The work by Luce and Raiffa (1958) is by now classical, as many examples traced back to this source, such as {\it Prisoners' Dilemma} and the {\it Battle of the Sexes}. However one of the first formal treatise on game theory is by Zermelo (1913), a logician who proved that in the game of Chess either White always wins or Black always wins or either a player can always force a draw.

\item Classical applications of game theory were in Economics, Behavioral sciences and Biology. Game Theory and its impact on the theory of Economics has led to several Nobel Prize winnings, starting with the 1072 Nobel Prize to Ken Arrow; followed by the 1994 to  John Nash, Reinhard Selter, John Harsanyi; then  Mirrlees and Vickrey (1996), Sen (1998), Akerlof, Spence and Stiglitz (2001), Aumann and Schelling (2005) Hurwicz, Maskin and Myerson (2007) and Roth and Shapley (2012). As a consequence, ideas, concepts and formal language of game theory form large parts of economics.

\item The most recent applications  are related to networked systems as reflected in Online advertisement on the Internet, Information evolution, Belief propagation in social networks; deployment of distributed passive and active sonars for underwater sensor field design. Game can be used indeed to model the interactions taking place in a network, with the network nodes, largely interdependent, acting as players competing or forming coalition to maximize their {\it quality of service.} For instance in the field of Signal Processing for communication networks one may design a game to address issues such as data security, spectrum sensing in cognitive radio, multimedia resource management and image segmentation. 
\end{itemize}

\end{remark}

\subsubsection{Existence and Computation of Equilibrium}

Nash (1951) famously proves that \textcolor{purple}{Every finite, noncooperative strategic game of two or more players has at least a (possibly mixed) Nash Equilibrium}. 

However, this is only an existence result; that is, it does not present the steps to find an equilibrium. The best response correspondence is one way for finding a NE.

For any $s_{-i}\in S_{-i}$ define the {\it best response correspondence} of player $P_i$ as 

\begin{equation}
B_i(s_{-i} := \{s_i\in S_i \quad |\quad u_i(s_i,s_{-i}) \ge u_i(s'_i,s_{-i})\quad \forall s'_i\in S_i\}
\end{equation}

A strategy profile $s^*$ is a NE if
$$
s^*_i \in B_i(s_{-i})\quad \forall i=1,\ldots,n.
$$
In other words, the Nash Equilibria are the fixed points of best-correspondence mapping
$$
B(s):= (B_i(s_{-i})_{i=1,\ldots,n}.
$$

\subsubsection{Mixed Strategy Nash Equilibrium}

Let consider a probability distribution $\Delta(S_i)$ over the pure strategy set $S_i.$ We call  $\sigma_i\in \Delta(S_i)$ the  mixed strategy of player $P_i.,$ and $\sigma=(\sigma_i)_{i=1,\ldots,n}$ is the mixed strategy profile. That is, $\sigma_(s_i)$ is the probability assigned to each strategy $s_i.$

The expected payoff/utility is defined as
\begin{equation}
 \tilde{u}_i(\sigma_i,\sigma_{-i}):=\mathbf{E}(u_i)= \sum_{s\in S} (\prod_{i=1,\ldots,n} \sigma_i(s_i) ) u_i(s)
\end{equation}

Therefore
\begin{definition}(Mixed NE)
A mixed strategy profile $\sigma^*$ is a \textcolor{blue}{Mixed Nash Equilibrium} (MNE) 
if
\begin{equation}
\tilde{u}_i(\sigma^*) \ge \tilde{u}_i(\sigma_i,\sigma^*_{-i})\quad \forall i=1,\ldots,n\quad \forall \sigma_i\in \Delta(S_i)
\end{equation}
\end{definition}

\begin{remark}(\textcolor{blue}{Pareto-Optimality})

The NE or any equilibrium is said to be {\it Pareto-inefficient} whenever there is another strategy profile in the game providing a payoff greater than the equilibrium payoff. This constitutes a drawback for considering a NE.

A strategy profile $s^*$ is therefore considered {\it Pareto-Optimum} if there is no other strategy profile $s$ such that $u_i(s)\ge u_i(s^*)\quad \forall i=1,\dots,n$
and $u_{i_0}(s) > u_{i_0} (s^*)$ for some $i_0.$
\end{remark}

That is to say that at a Pareto-Optimum NE, the payoff to one player cannot be increased without decreasing the payoff to at least one other player. A classic example of an inefficient NE is in the well-known {\it Prisoner's Dilemma}: The unique NE is {\it Pareto inferior}: indeed, if both players choose (cooperate, cooperate), they would both get a higher payoff.

In addition to these characteristics an equilibrium solution may well feature other socially sounding characteristics, such as a {\it maximizer of social welfare,} originated in the application of game theory to economics. (Debreu 1952). Indeed

\begin{definition}
A strategy profile $s^*$ is called \textcolor{blue}{social optimum} if it satisfies
$s^* \in \arg \max _{s\in S} \sum_{i=1,\ldots,n} u_i(s) $
\end{definition}

As one of the most important solution concepts, the efficiency of the Nash Equilibrium has to be constantly assessed. It can also be measured using the following concept of \textcolor{blue}{Price of Anarchy} (PoA) defined as

\begin{equation}
PoA:= \frac{\max _{s\in S} \sum_{i=1,\ldots,n} u_i(s)}{\min _{s\in S^*} \sum_{i=1,\ldots,n} u_i(s)},  
\end{equation}
where $S^*$ is the set of all Nash Equilibria in the game.

That is to say that the {\it Price of Anarchy} measures the performance loss of the worse NE when compared to a socially optimum strategy. The efficiency of the NE is higher whenever the PoA is closer to 1. 

{\it Mechanism Design} refers to a process of improving the NE efficiency by transforming the game while keeping the NE, noting, for instance, that payoffs are invariant under {\it affine transformation.} Of course, one may also improve efficiency by modifying the solution concept with the game unchanged, Such modification may refer to equilibrium concept such as \textcolor{blue}{Correlated Equilibrium} (CE) or the \textcolor{blue}{Nash Bargaining Solution}(NBS), both here recalled, See also Nash 1951 and Aumann 1959 in the references.

\begin{definition}(\textcolor{blue}{Correlated Equilibrium})

A {\it Correlated Equilibrium} is a joint probability distribution $\sigma^*\in \Delta(S)$ satisfying
\begin{equation}
 \sum_{s_{-i}\in S_{-i}} \sigma^*(s_i,s_{i}) u_i(s_i,s_{-i}) \ge  \sum_{s_{-i}\in S_{-i}} \sigma^*(s'_i,s_{i}) u_i(s'_i,s_{-i}),  
\end{equation}
 where $i=1,\ldots,n$ and $ s'_i\in S_i.$   
\end{definition}
See more details in Aumann's work.(Aumann 1959): A {Traffic light} with its 3 distinct alternating colors of red, yellow and green, acts as a trusted neutral third party coordinating device, is considered as a simple example of correlated equilibrium: no driver has an incentive to disobey a traffic light while everyone else is following it. It serves as a joint probability distribution correlating the actions of all drivers, with a high probability of one driver seeing green while the other driver sees red.

The CE is said to be {\it Coarse} if it requires an {a-priori} commitment to follow the given recommendation, e.g.,  getting a driving license implies an a-priori commitment to obey a traffic  light while driving. 
We indicate below how {\it entanglement} of quantum strategies could be considered strong non-classical correlation between strategies.

Seeking an alternate efficient may also lead to the so-called {\it Nash Bargaining Solution}(NBS) defined here for a two-player game as

\begin{definition}(\textcolor{blue}{Nash Bargaining Solution}) (NBS)

The Nash Bargaining Solution is the unique Pareto-Optimal solution to
\begin{align}
&\max_{(u_1,u_2)\in \cal{U|}}(u_1-\lambda_1)(u_2-\lambda_2)\\
&\mbox{subject to}\quad  u_1 \ge \lambda_1,\quad u_2 \ge \lambda_2,
\end{align}
where $\cal{U}$ is a closed convex set of utility points. The product $(u_1-\lambda_1)(u_2-\lambda_2)$ is called the {\it Nash product.}
\end{definition}

\subsubsection{ Potential and Congestion games}

Many game scenarios are concerned with finding a {\it pure} Nash equilibrium. There are two well-known classes of games which always possess at least one pure Nash equilibrium. We recall the following result in the case of infinite strategy set games,

\begin{theorem}(Debreu, Fan, Glicksberg)
Consider a strategic form game with an infinite strategy set such that
\begin{enumerate}
\item $S_i$ is convex and compact
\item $u_i(s_i,s_{-i})$ is continuous in $s_{-i},$ continuous and quasiconcave in $s_i.$
\end{enumerate}
Then the game has a pure Nash Equilibrium. 
\end{theorem}
See proof and more details in Debreu, Fan or Glickberg 1952

\begin{definition}(Potential Games)

A game is a {\it potential game} if there exists a potential function $\Phi,$ real-valued over the strategy set $\mathbb S$ such that for every player $P_i$ and for every strategy profile $s=(s_i,s_{-i})$ and every strategy $s'_i$ we have 
\begin{equation}
    \Phi(s_i,s_{-i}) - \Phi(s'_i,s_{-i}) = u_i(s_i,s_{-i}) - (s'_i,s_{-i})
\end{equation}
\end{definition}
That is, in potential games the change in the player's utility upon switching of strategy remains the same as the change in the potential function. (Rosenthal 1973).

It was then proved that

\begin{theorem}
Every potential game has at least one pure Nash equilibrium, namely the strategy that maximizes the potential function $\Phi.$
\end{theorem}

Here also there is no indication on how to determine that pure equilibrium. However, the proof of this theorem reveals that the {\it best response dynamics} leads to finding the pure Nash equilibrium, without saying much about the number of steps. That is
\begin{theorem}
In any finite potential game, best-response dynamics always converge to a pure Nash equilibrium.
\end{theorem}

\begin{definition}(Congestion Games)

A {\it congestion game} is defined by a finite set $P$ of n players, a finite set $S$ of resources. Each player $P_i$ has a set $S_i\subset S$ of allowable strategies. For each resource $r_i\in S$ and assuming $x_j$ players using $r_j$, there is a load dependent cost function $c_j(x_j)$ so that the total cost to a player $P_i$ is $\sum_{j\in S_i} c_j(x_j)$.
\end{definition}

The simplest example of congestion game is the non-atomic network routing game. It was proved that
\begin{theorem}(Rosenthal 1973)
Every congestion game is a potential game
\end{theorem}
The implication is, therefore, that every congestion game has a pure Nash found using best response dynamics.

\subsubsection{Mean Field Game}

The study of decision-making in large population of interacting agents is much more involved than in the above scenarios, and it is the subject of the sub-game theory of \textcolor{blue}{Mean Field Game} (MFG), introduced by Jean-Michel Lasry and Pierre-Louis Lions and others in the 2000s. (Lions and Lasry, 2007). Its focus is on the collective or average behavior of the game population rather than the individual one.

The {\it Hamilton-Jacobi-Bellam equation} and the {\it Fokker-Planck equations} are some of the most common equations studied in MFG, respectively addressing the optimal control of every single player and the evolution over time of the distribution of players. (Djehiche et al., 2017). These fundamental equations of MFG have led to some interesting insights on existence and uniqueness conditions for equilibria.  They are respectively

$$
-\partial_t u(x,t)+ H(x,\nabla u) = f(x,m)\quad \text{(\textcolor{purple}{Hamilton-Jacobi-Bellam})}
$$

where $u(x,t)$ is the value function, $H$ the Hamiltonian and $f(x,m)$ the cost function in terms of population density $m(x,t).$ It is considered a "backward PDE".

$$
\partial_t m - \nabla \times (m\nabla H_p(x,\nabla u)) =0,\quad \text{(\textcolor{purple}{(Fokker-Planck)}}
$$

with $H_p$ the gradient of the Hamiltonian with respect to momentum accounting for the drift in the players dynamics. It is considered a "forward PDE" on the population distribution.

Solutions methods are varied, and oftentimes ad hoc, and include: fixed-point techniques. Variational schemes (convex optimization), and of course, Numerical schemes (Finite Difference, Finite Elements, neural-network-based solvers).

The following Fixed Point Theorems are widely used in MFG:

\begin{theorem} (\textcolor{blue}{Schauder Fixed Point Theorem})(Schauder, 1930)

If $\mathbb X$ is a Banach space and $T:\mathbb X \rightarrow \mathbb X$ is a compact, continuous operator mapping a convex, closed, and bounded subset of $\mathbb X$ into itself, then $T$ has a fixed point.
\end{theorem}
and also

\begin{theorem} (\textcolor{blue}{Banach Fixed Point Theorem}) (Banach, 1922)

If $(\mathbb X,d)$ is a complete metric space and $T:\mathbb X \rightarrow \mathbb X$ is a contraction mapping (i.e., there exists $0 \le c < 1$  such that $d(T(x),T(y)) \le c d(x,y) $ for all $x,y \in \mathbb X$), then $T$ has a unique fixed point $x^*$ such that $T(x^*) = x^*.$

\end{theorem}

In some notable recent developments, MFG provides an ideal framework for modeling financial markets and resource allocation, traffic flow and communication networks, epidemics, and multi-agent reinforcement learning and  swarm intelligence. (Yves, 2020 and Carmona, 2020). 

MFG has been experiencin growing relevance in analyzing and proposing equilibrium solutions for complex systems with large numbers of interactive decision-makers, leveraging powerful computational algorithms. It has been extended to the so-called Ergodic MFG to investigate long-term behavior in large population interactive decision-making, with time-independent solution equilibrium.

\subsubsection{Approximate Nash Equilibrium}

Recall that Nash's famous theorem is only an existence theorem about Nash Equilibrium (NE), based on Brouwer and Kakutani's fixed point theorems, respectively recalled as

\begin{theorem} (\textcolor{blue}{Kakutani Fixed Point Theorem}) (Kakutani, 1941)

Let $\mathbb K \subset \mathbb{R}^n$ be a compact and convex subset and $F:\mathbb K\rightarrow 2^{\mathbb K}$ a mapping into the set of compact and convex subsets of $\mathbb K$ which is \textcolor{yellow}{upper semi-continuous}\footnote{$F$ is upper semi-continuous in $x\in \mathbb K,$ if $\forall (x_k)_{k\in \mathbb N} \in \mathbb K\quad x_k\rightarrow x $ and $\forall (y_k)_{k\in\mathbb N} \in \mathbb K,\quad y_k\in F(x_k), \quad y_k \rightarrow y$ then $y\in F(x)$. } in every $x\in\mathbb K,$ then $\mathbb F$ possesses a fixed point, i.e., there is some $\hat{x}\in \mathbb K$ with $\hat{x}\in F(\hat{x}).$
\end{theorem}

It is a generalization of the \textcolor{blue}{Brouwer Fixed Point Theorem} stating:
\begin{theorem}(Brouwer, 1911)

Every continuous mapping of a convex and compact subset of an n-dimensional Euclidean space into itself.
\end{theorem}

The computation of the NE  remains challenging and complex even with the advances in computational methods assisted by AI/ML. To date no efficient algorithm is known. Known algorithms include:

\begin{itemize}
    
\item {\it Regret Minimization}\footnote{Popularized by Jeff Bezos after leaving his well-paid hedge fund job to start Amazon, as he might have regretted not pursuing his entrepreneurial vision.} (no-regret learning): strategies are adjusted taking into account players past performance to minimize regret

\item {\it Reinforcement Learning}: Decision-making through interactions with environment to include merit-based rewards/penalties to maximize cumulative rewards over time.

\item {\it Linear Programming,} in particular, for zero-sum games, in an attempt to compute near-optimal strategies: determine optimal solution to problems with linear relationships, constraints and objective function (e.g., maximize profit and minimize cost)

\end{itemize}
As it is often the case in mathematics, one way around the problem has been the notion of {\it Approximate Nash Equilibrium} (ANE) to some $\epsilon$ accuracy. That is, for a game in a pure or mixed normal form, one would like to evaluate a minimal tolerance value to express how far a strategy profile $\hat{s}$ is from the exact NE $s^*$ whose existence is guaranteed in the Nash Existence Theorem. In other words, 

\begin{definition}
 A strategy profile $\hat{s}=(\hat{s}_i,\hat{s}_{=i}) $ is an $\epsilon-$Nash Equilibrium ($\epsilon \ge 0)$ if no unilateral deviation for the exact NE $s^*$ leads to an expected payoff gain of  more that $\epsilon$. Mathematically,

 \begin{equation}
u_i(\hat{s}) \ge u_i(s_i,\hat{s}_{-i}) - \epsilon,\quad \forall i=1,\ldots,n\quad \forall s_i\in S_i
 \end{equation}
 
\end{definition}

 Due to the intractability of computing Nash equilibria, approximations seem more practical and feasible,  as it is in general in numerical analysis or approximation theory. Moreover, in real-world scenarios, games of chance or games of strategy, agents/players act optimally within a "reasonable" margin error, despise the "rationality" assumption.

We now turn to one of the most important class of games.

\section{Evolutionary Game Theory}

\subsection{Introduction and Tenets of EGT}

Evolutionary Game Theory (EGT)\footnote{EGT first appeared in the work "The Logic of Animal Conflict" of Maynard Smith, a mathematical biologist and George Smith, attempting to understand ritualized animal behavior.} was adapted from the biological evolution theory and developed the mathematical theory of Darwinian evolution, providing the most promising applications of the theory of games: it has a profound impact on social and behavioral sciences from economics to psychology to linguistics and language evolution to learning and cooperation in agent-based systems (origin of social norms and cultural trends). EGT contributes to addressing questions such as whether a fly, a fig tree or the Harper experiment ducks, all be decision-makers evaluating all possible outcomes to select their optimal strategy for "success".

From the early works of Smith and Price statements, the subsequent theory of evolution was based on {\it games of selfish genes} rather than considering organisms acting {\it for the good of the species}. (Smith, 1982)

We complete this introduction of EGT with some quotes from Malthus and Darwin
T.R. Malthus:{\it Through the animal and vegetable kingdom, nature has scattered the seeds of life abroad with the most profuse and liberal hand; but has been comparatively sparing in the room and nourishment necessary to rear them.}

To which, Darwin reacted: {\it Fifteen months after I had begun my systematic enquiry, I happened to read for amusement "Malthus on Population"$\ldots$ it at once struck me that $\ldots$ favorable variation would tend to be preserved and unfavorable ones to be destroyed. Here, then, I had at last got a theory by which to work}

\begin{remark}

EGT mainly considers a large population of individuals, the players in the broader sense. These players $P_i$ have a finite number $N$ of strategies $s_{i=1,\ldots,n}$ each of which has a frequency or a probability of occurrence that changes over time in response  to the decisions of all the other players under the principle that individuals switch to those strategies with better payoffs from poorly performing strategies. 

The changes in frequency of the strategies in the population are described by a probability distribution $\Delta(S)$ on the finite set $S=(s_{i=1,\ldots,n})$ of the strategies. 
\end{remark}

\subsection{Evolution Matrix Game}

Consider a large population whose individuals have $n$ strategies $s_i\in S$ for interaction with each other, together with a probability distribution $\Delta(S)$ over the finite  set of strategies. That is, let $p_i\in[0,1]\subset \mathbf R,$ be the probability for strategic value $s_i$ to be chosen. The vector of probabilities or frequencies $p=(p_1,\ldots,s_n)\in [0,1]^n\subset \mathbf R^n$ describes a state of the population as well as a behavioral strategy profile , with the components $s_i$ satisfying the probability constraint $\sum_i p_i=1.$ The state vectors $p=(p_{i=1,\ldots,n}$ lies on the \textcolor{blue}{(n-1) simplex space} of strategies

\textcolor{purple}{\begin{equation}
\Delta=\{s=(p_1,\ldots,p_n)\in \mathbf R^n\quad |\quad 0\le p_i\le 1,\quad \sum_i p_i=1\}
\end{equation}}

compact and convex spanned by the set of vectors $e_i$ of the standard basis of the Euclidean space $\mathbf R^n.$
\begin{equation}
e_i=(0,\ldots,0,1_i,0,\ldots,0),\quad i=1,\ldots,n.
\end{equation}

\begin{remark}
\noindent
\begin{itemize}

  \item Pure populations states are given by the corners vectors of the standard Euclidean basis $e_i=(0,\ldots,0,1_i,0,\ldots,0),i=1,\ldots,n.$ These states are also called {\it homogeneous}. 
  
  \item Note that the probability $p_i$ assigned to strategy $s_i$ changes with respect to payoffs by learning, copying, inheriting etc. And in turn payoffs are function of the probabilities yielding a {\it feedback loop dynamic}. The difference $h_i(s)$ has been identified in certain "game dynamics" with a relative increase of the frequency.
  
  \item Strategies with higher payoffs reproduce faster. Consequently individuals have a tendency to vary the probability assigned to a core value based on some perceived relative advantage.
  
  \item Subscribing to competing strategic values lead to the emergence of {mutant subgroups}, with the by-products of the appearance of a variety of life forms, beliefs, cultures, languages, practices or techniques.
  
  \item In contrast to other types of dynamics such as Replicator, Replicator-Mutator, or Logit Dynamics, any kind of revision protocol is described by changes in the values of the probability assigned to each strategy, not by switching of strategy. Successful strategies are represented by "successful individuals" and could be learned (imitated), possibly inherited. (Hofbauer et al., 1979, 1998, 2003).
  
  \end{itemize}
\end{remark}

The interaction between an $s_i$-individual and an $s_j$-individual results in a {\it payoff}/{\it reward}/{\it utility} denoted by $a_{ij}\in\mathbf R.$ Therefore all the payoffs form the {\it evolution payoff matrix}

\begin{equation}
A=(a_{ij})_{1\le i,j\le n}
\end{equation}

defining the game as an evolution matrix game. Based on the context, the real entries of this $n \times n $ matrix could be assumed to satisfy various properties of matrices. That is, the matrix $A$ could be {\it symmetric}, i.e., $A=A^T;$ {\it skew-symmetric}, $A=-A^T;$ {\it cyclic symmetric} (Toeplitz circulant with indices counted cyclically modulo n); or a {\it banded matrix} (sparse with non-zero entries confined to a diagonal band).
These features may help to further analyze the matrix

Accordingly we define a payoff or utility function $u=(u_1,\ldots,u_n):\Delta\subset [0,1]^n\rightarrow \mathbf R^n,$ also called the {\it payoff vector field}. 

\begin{definition}
\noindent
\begin{enumerate}
  \item Given a population state determined by the probability vector $p=(p_1,\ldots,p_n)\in\Delta,$ the {\it expected payoff} of an $s_i$-player is defined as
 
 \textcolor{purple}{\begin{equation}
u_i(p)=e_i A p^T=\sum_{j=1}^n a_{ij} p_j=(A p)_i.
\end{equation}}

The payoff functions $u_i$ are continuous linear functions. 
  \item The {\it average payoff/utility} in the population state $p$ is given by

\textcolor{purple}{\begin{equation}
\bar{u}(p)=p A p^T=\sum_{1\le i,j\le n}p_i A p_j.
\end{equation}}

  \item We denote by 
\begin{equation}
h_i(p)=u_i(p) - \bar{u}(p),
\end{equation}

that is, the excess payoff between the individual expected payoff and the average payoff in the population state, which actually impacts the frequency and the probability value of a given strategy.

\end{enumerate}
\end{definition}

Note that the sign of excess payoff $h_i(p)$ for the strategy $s_i$ in the population state $p$ determine its variation in frequency as well as its corresponding probability. In the {\it standard or regular} population state $p$ each strategy $s_i$ receives a nonzero probability, that is, $p_i\in (0,1).$ We denote for a general $n$ the support of a state $p$ by 
$$
supp(p)=\{i\in\{1,\ldots,n\}|\quad p_i>0\}
$$
We now define the equilibrium of such an evolutionary game. 

\begin{definition}
A population state $p^*=(p^*_1,\ldots,p^*_n)\in\Delta$ is called a {\it Nash Equilibrium,} denoted NE, if

\textcolor{purple}{\begin{equation}
\bar{u}(p^*)=p^* A p^{*T}\ge p A p^{*T},\quad \forall p\in\Delta.
\end{equation}}

Equivalently we have
$$
\aligned
\bar{u}(p^*)&=p^* A p^*=max_{i=1,\ldots,n} u_i(p^*)\\
u_i(p^*)&=u_i(p^*_i,p^*_{-i})\ge u_(p_i,p^*_{-i})
\endaligned
$$
where again 
$$
u_i(p_i,p^*_{-i})=u_i(p^*_1,\ldots,p^*_{i-1},p_i,p^*_{i+1},\ldots,p^*_n),\quad \forall p_i\in\Delta_i=[0,1].
$$
\end{definition}

\begin{remark}[Interpretation]

In the Nash Equilibrium state $p^*$ the unilateral change of the probability of a single strategy does not lead to a higher payoff. Under the rationality assumption of the game one would expect such a state to be maintained. That is not in general the case because individuals do not always behave rationally. Deviations are expected dictated by changes in preferences, and the stability of the Nash Equilibrium state $p^*$ is not guaranteed. We are therefore interested in dynamic conditions that could ensure the {\it evolutionary stability} of that equilibrium as defined next.
\end{remark}

\begin{definition}[\textcolor{blue}{Evolutionary Stability}]
A Nash Equilibrium state $p^*\in\Delta$ is {\it evolutionarily stable} if
\begin{equation}
\bar{u}(p^*)=p^* A p^* = p_0 A p^{*T}\quad\text{for some $p_0\in\Delta,$ $p_0\ne p^*$}
\end{equation}
implies
\begin{equation}
\bar{u}_(p_0)=p_0 A p_0 < p^* A p_0^T.
\end{equation}
\end{definition}

That is, if the same payoff of the NE state $p^*$ could be achieved in some other state $p_0,$ then this other state cannot be a NE state. Conditions are needed to characterize dynamically both states, necessary and/or sufficient conditions. Here are some results widely proved in the literature on Evolutionary Game Theory. See for instance Hofbauer and Sigmund 1998 or Hofbauer and Sigmund 2003.

\begin{itemize}

\item Let $p^*\in\Delta$ be a NE state with support $supp(p^*)$. Then we get
$$
u_i(p^*)=\bar{u}(p^*)\quad \forall i\in supp(p^*)
$$
In a NE state, the individual expected payoff is the same as the average in the state for the nonzero probabilities making the NE state.

\item As indicated above the payoff functions $u_i$ are continuous linear functions on the (n-1)-simplex $\Delta,$ a convex compact subset of $\mathbf R^{n}.$ This ensures 
$\forall p^*\in\Delta,$ $\exists \tilde{p}_i\in\Delta_i=[0,1],$ such that 
$$
u_i(\tilde{p}_i,p^*_{-i})\ge u_i(p_i,p^*_{-i})\quad | \quad \forall p_i\in \Delta_i
$$
\end{itemize}

\begin{remark}
The population game as described above is so far in a {\it static form.} Players  as rational decision-makers adapt their strategies over time by varying the probabilities assigned to each strategy  for a more efficient distribution of their preferences and behaviors. Indeed strategies with higher payoff/rewards spread quickly in the population through, e.g.,  learning, copying/imitating. Payoffs depend on the frequencies/probabilities of strategy in the population, which themselves change according to the payoffs/rewards, resulting in a feedback loop dynamic characterizing the evolution of the population. We therefore introduce a dynamic in the game to account for the changes over time of the probability distribution $p(t)=(p_i(t))_{i=1,\ldots,n}$
\end{remark}

 \subsubsection{Towards the Replicator Dynamics}

The per capita growth rate of the time-continuous dependent probability $p_i(t)$ is the logarithm derivative 
\begin{equation}
\frac{d(log(p_i(t)))}{dt}=\frac{\dot p_i}{p_i}
\end{equation}

is determined by the difference between the expected payoff to the core value $A_i$ and the average payoff to the population state, leading to the time-continuous dynamical system or \textcolor{blue}{Replicator Dynamics} $(\cal{RD})$
\begin{equation}
\frac{\dot p_i(t)}{p_i}=h_i(p)=u_i(p)-\bar{u}(p),
\end{equation}

or equivalently
\textcolor{purple}{\begin{equation}
\dot p_i(t)=p_i(u_i(p)-\bar{u}(p))=p_i((Ap)_i-p A P^T)=g_i(p),
\end{equation}}

\begin{remark}
\noindent
\begin{enumerate}

  \item This is a system of differential equations for $i=1,\ldots,n$ It appears in various areas such as population genetics, chemical networks, and is famously known as the {\it replicator dynamics}\footnote{A replicator is an entity with the means to make accurate copies of itself. It can be a gene, an organism, a belief, a technique, a convention, a cultural norm.} introduced by Taylor and Jonker (1978) and coined as such by Schuster and Sigmund (1983).
  
\item Under the constraint $\sum_{i=1}^n p_i(t) =1$ the system of equations actually reduces to (n-1) differential equations we analyze on a (n-1)-simplex. 

\item Moreover adding a constant to a column entries of the payoff matrix $A$ does not change the equations and its dynamic properties. One may then set the diagonal entries to zeros, or set the last row of $A$ to zero in analyzing the dynamics. Such simplifications are indeed undertaken in practice.
\item The well-defined power product $V(p)=\prod_i p_i^{\alpha_i}$ satisfies 
\begin{equation}
\dot V = V\sum \alpha_i[(A p)_i - p A p^T].
\end{equation}

\item The frequency of a pure strategy in the society as given by the probability increases when it has above average utility/payoff. Individuals have limited and localized knowledge of the whole system, according to the {\it distribution of information}. Some strategies or core values could become extinct as time goes to infinity. Indeed whenever a core value $s_i$ is recursively strongly dominated as defined above, it will not survive the evolution.  That is, for the assigned probability $p_i(t)$, we have $\lim_{t\rightarrow\infty}p_i(t)=0.$ The corresponding core value therefore goes extinct. This process might explain the steady variation of strategic values and norms.

\end{enumerate}

\end{remark}
Fundamentally it is shown, e.g., in Weibull, 1995 that 
the mean average payoff or utility of a population state increases along any trajectory of the evolution equation $(\cal{RD})$ giving by the equation
\begin{equation}
\dot{\bar{u}}=2\sum_{i=1}^n p_i (u_i(p)-\bar{u})^2.
\end{equation}

It is similar to a result by the biologist Fisher with respect to Natural Selection. (Fisher 1922)

We analyze the system on the unit (n-1)-simplex subset of $\mathbf R^{n},$ and the state space of the vectors of probabilities. It is actually a cubic polynomial dynamic, whose class has been extensively investigated. We denoted $\dot{\Delta}=int(\Delta)$ the {\it interior} of the simplex given by the set
\begin{equation}
\dot\Delta :=int(\Delta)=\{p=(p_1,\ldots,p_n)|\quad p_i>0\quad \forall i=1,\ldots,n\}.
\end{equation}

The {\it boundary faces} are denoted by $\partial \Delta$ and defined as
\begin{equation}
\partial \Delta(J)=\{p\in\Delta: p_i=0,\quad \forall i\in J\},
\end{equation}
where $J$ is any nontrivial subset of $\{1,\ldots,n\}.$
Note that the hyperplanes $\sum p_i=1$ and $x_i=0$ are invariant, as well as the simplex faces. See Weibull 1995, Hofbauer et al. 2003.

The rest points are the zeros of $g(p)=(g_i(p))=0.$ That is, points $p\in\Delta$ such that $(Ap)_i=p A p^T,$ for all $i\in supp(x).$ Thus an interior rest point is a solution of the system of linear equations 
\begin{equation}
(Ap)_1=(Ap)_2=\ldots=(Ap)_n.
\end{equation}
The following so-called \textcolor{blue}{Folk Theorem of Evolutionary Game Theory} contains all relevant results of the Replicator Dynamics. 

\begin{theorem}
\noindent
\begin{enumerate}
  \item If $p^*$ is a Nash Equilibrium then it is a rest point for $(\cal{RD})$
  
  \item If $p^*$ is an evolutionarily stable NE then it is asymptotically stable for $(\cal{RD})$ globally if $p^*$ is interior.
  
  \item If the rest point $\hat p$ is stable then it is a Nash Equilibrium
  
  \item Every interior rest point $\hat p\in \dot \Delta$ is a Nash Equilibrium 
  
  \item If a rest point $\hat p$ for $(\cal{RD})$ is also the forward limit point of an interior orbit $x(t)$ of $(\cal{RD}),$ i.e., $x(t)\in \dot\Delta=int(\Delta)$ then $\hat p$ is a Nash Equilibrium.
  \end{enumerate}
\end{theorem}

\begin{remark}

The converse is not true in the {\it Folk Theorem of Evolutionary Game Theory.}
\end{remark}

Additional important well-known results (Hofbauer et al. 2003)  of EGT are summarized as follows.

\begin{remark}
For a boundary rest point $\hat p$ the difference $h_i(\hat p)=(A\hat p)_i-\hat p A \hat{p}^T$ is also an eigenvalue $\lambda_{\hat p}$ for the Jacobian $J(\hat p)$ with an eigenvector transversal to the face $\hat p_i=0.$ That entails
A rest point $\hat p$ for the evolution equation $(\cal{RD})$is a Nash Equilibrium if and only if its transversal eigenvalues are non-positive.
\end{remark}

Another important question is the following: Are these strategies {\it permanent and resistant} to shocks of any nature and size, random or otherwise? It amounts to the notion of {\it permanence} or {\it self-sustained} for the {\it replicator dynamics} in the sense:
\begin{definition}
The evolution equation $(\cal A)$ is {\it self-sustained} if $\exists$ a compact set $K\subset int(\Delta)$ such that $\forall x\in int(\Delta)$ $\exists$ a real time $T$ such that $p(t)\in K\quad \forall t>T.$
\end{definition}
 It then results
 
\begin{theorem}

If the evolution equation $(\cal{RD})$ is {\it self-sustained} then it has a unique rest $p^*$ in the interior $\dot\Delta$ which is therefore the unique national Nash Equilibrium $NE_A$.

Moreover along each interior orbit $p(t)$ the {\it time averages} converge to $p^*$ 
\begin{equation}
\bar{a}_i:=\frac{1}{T}\int_0^T p_i(t) dt \longrightarrow p^*_i,
\end{equation}

as $T\rightarrow \infty,$ and $i=1,\ldots,n.$
\end{theorem}
 The system is said to be {\it ergodic.} See again Hofbauer et al. 2003
 
 Consequently, if an interior orbit $p(t)$ has a boundary $\omega-$limit point the time average need not converge. And every orbit converges to the boundary 
$\partial\Delta$ in the absence of an interior rest point.

It will be very interesting to analyze the dynamics of these subgames (i.e., using a reduced number of pure strategies) and their relation to the overall dynamics, with the possibility of appearance of special dynamics such as {\it limit cycles.} 

\begin{remark}
Before presenting in the following sections on the most recent advances in game theory, we summarize here the most salient characteristics of the Classical Game Theory
\end{remark}

\subsection{Principles of Classical Game Theory}

We recall first the three most common models of classical game theory, widely described and studied in most books on game theory:

\begin{enumerate}
\item \textcolor{blue}{The Prisoner's Dilemma} (PD): This is a paradigm of noncooperative game, which mathematically models how players' mutual cooperation with a better collective outcome/payoff is not a Nash Equilibrium of the game. But instead, mutual betrayal/defection leads to the unique Nash Equilibrium, which is not event Pareto-optimal. That is, the tension between individual rationality and collective benefits. Its applications include Price Wars and Trade Agreements, Cartels operations, Arms races, Plea Bargains, Machine Learning and Artificial Intelligence Ethics. The PD has the following general payoff matrix

\begin{table}[h!]
    \setlength{\extrarowheight}{2pt}
    \begin{tabular}{cc|c|c|}
      & \multicolumn{1}{c}{} & \multicolumn{2}{c}{Column Prisoner}\\
      & \multicolumn{1}{c}{} & \multicolumn{1}{c}{$C$}  & \multicolumn{1}{c}{$D$} \\\cline{3-4}
      \multirow{2}*{Row Prisoner}  & $C$ & $(R,R)$ & $(S,T)$ \\\cline{3-4}
      & $D$ & $(T,S)$ & $(P,P)$ \\\cline{3-4}
    \end{tabular}
  \end{table}
with $T>R>P>S,$ and classically, T is the temptation to defect, R the reward for mutual cooperation, P the punishment for defection, and S the so-called "sucker"'s payoff to the player who defects while the other cooperates.

\item \textcolor{blue}{Matching Pennies}: This models zero-sum games with no pure Nash Equilibrium, with the typicall payoff
\begin{table}[h!]
    \setlength{\extrarowheight}{2pt}
    \begin{tabular}{cc|c|c|}
      & \multicolumn{1}{c}{} & \multicolumn{2}{c}{Column Player}\\
      & \multicolumn{1}{c}{} & \multicolumn{1}{c}{$H$}  & \multicolumn{1}{c}{$T$} \\\cline{3-4}
      \multirow{2}*{Row Player}  & $H$ & $(1,-1)$ & $(-1,1)$ \\\cline{3-4}
      & $T$ & $(-1,1)$ & $(-1,-1)$ \\\cline{3-4}
    \end{tabular}
  \end{table}

  The optimal strategy is mixed, i.e., randomization with equal probability, leading a mixed Nash Equilibrium. Its applications include Markets Competition, Bidding Wars, Hacker Vs. Defender, Submarine warfare, and Soccer Penalty Kicks.
  
\item \textcolor{blue}{The Battle of the Sexes}: A paradigm for coordination games, with a payoff matrix such as
\begin{table}[h!]
    \setlength{\extrarowheight}{2pt}
    \begin{tabular}{cc|c|c|}
      & \multicolumn{1}{c}{} & \multicolumn{2}{c}{Husband}\\
      & \multicolumn{1}{c}{} & \multicolumn{1}{c}{$O$}  & \multicolumn{1}{c}{$F$} \\\cline{3-4}
      \multirow{2}*{Wife}  & $O$ & $(3,2)$ & $(0,0)$ \\\cline{3-4}
      & $F$ & $(0,0)$ & $(2,3)$ \\\cline{3-4}
    \end{tabular}
  \end{table}
  Players with conflicting individual preferences benefits from being together. It has two pure Nash Equilibria and one mixed Nash Equilibrium. It can be applied to Communication and Negotiation scenarios, Trade Agreements Military Alliances, Business and Technology Standards,

\end{enumerate}

\begin{remark}
   As presented in the subsequent sections, these games will greatly benefit from  straight \textcolor{blue}{quantumization and p-adic quantumization} as described below.   
  \end{remark}
  
Classical game theory, in particular, in its multi-player strategic form, mathematically describes the strategic interaction, the {\it game}, conflictual or cooperative, between two or more entities, the {\it players}. A player's {\it pure strategies} refer to their ability to interact with the other players, and their stakes in the game are their {\it payoffs or utilities,} which, as {\it rational entity}, they seek to optimize/maximize, consistent with some preferences ordering over their payoffs.

A {\it game play} calls for the choice of a pure strategy by each player, leading to a tuple of pure strategies, the {\it pure strategy profile}, which determines the player's payoff, under the assumption that each player's choice of pure strategy is indeed their {best response} to the others' choice. The outcome of choosing and using such pure best response strategy is the pure strategy {\it Nash Equilibrium} (PNE), a fundamental goal of multi-player game theory.

Importantly, the appearance, if any, (as it may not even exist) of a PNE does not imply {\it uniqueness or optimality}. Players may (and in general, they must) therefore extend their strategy set to include {\it mixed strategies,} (von Neumann) a randomization between their pure strategies through an efficient probability distribution over the pure strategy set, and yielding {\it expected payoffs.} 

A fundamental result for the strategic form (with a finite strategy set) multi-player game is Nash's famous theorem, we recall as 

\begin{theorem}(Nash Existence Theorem)
If each player in a multiplayer game (two or more players) has a finite number of pure strategies, then the game has a Nash Equilibrium (not necessarily unique) in (possibly) mixed strategies.
\end{theorem}

The extension of a game strategy or payoff profile can take several possible directions. In the classical case, the co-domain of the payoff function is the real numbers field, $\mathbb{R}$, which can be switched to the p-adic numbers field $\mathbb{Q}_p$ described below. Another possible extension, exposed below, is the use of \textcolor{blue}{quantum strategies}, initiated by Meyer in 1999.

Additional extensions are also presented below in the sections on \textcolor{blue}{Intelligent Games} and \textcolor{blue}{Non-Archimedean or p-adic Games}.
     
\section{Intelligent Game Theory}

In addition to the above categories of games, the so-called \textcolor{blue}{Intelligent Game Theory} (IGT) is a fairly recent extension. It combines game theory with the advancement of machine learning (ML), artificial intelligence (AI) and computational methods to improve decision-making in strategic environments, to learn, adapt and optimize accordingly. This leads to applications in robotics, autonomous systems, multi-agents systems and cybersecurity and many other emerging areas.

\subsection{ Multi-Agents Systems and Learning in Game theory}

Constructing agents that exhibit intelligent behavior is an expanding endeavor in research on computer systems, artificial intelligence, agent-based software development. So intelligent agents  refer to smart programs able to learn from their environments and user interactions to perform tasks and make their own decisions in order to maximize overall benefits (utility-based agents), to achieve a common goal (multi agents in transportation systems, robotics, social networks, ...).

Intelligent Game Theory aims at designing efficient strategies for such agents to leverage their own decisions against the actions of other agents in a dynamic and potentially adversarial environment.

Through Reinforcement Learning, these agents learn optimal strategies using trials and errors, rewards and penalty in order to continuously adapt their strategies and optimize outcomes. This is especially useful in competitive games, traffic management, network optimization.

To handle high-dimensional state spaces and higher order complexity the intelligent agents approximate complex game-theoretic models using advanced ML techniques such as deep neural networks.

It goes without saying that intelligent agents operate under conditions of uncertainty or incomplete information. Bayesian game  theory, approximate inference methods, belief-based learning assist these agents  to optimize tier strategies in such environments. AI-driven strategies can help predict and optimize outcomes including in iterated, repeated and incomplete games,

This is to say that there are ample and varied opportunities to apply Intelligent Game Theory. We name here the following:

\begin{itemize}

\item Social Networks: we now live in the age of expanding social networks on platforms such as YouTube, Instagram, Facebook, X/Twitter, Tik-Tok. AI assisted Game Theory is used to model interactions such as the so-called {\it reputation systems} where the actions and behaviors of others influence and impact sometimes negatively most users' behavior. Combining game-theoretic models with AI assisted recommendation systems and sentiment analysis will help predict and influence user behavior.
\item Healthcare: This is one important area of application of Intelligent Game Theory. Decision-making systems involve a variety of stakeholders: patients, doctors and their medical assistants, medical insurers. Here intelligent agents learn optimal strategies for resource allocation in medical treatments. This will lead to a complete and positive overhaul of the healthcare system, in particular, in countries with no universal efficient healthcare system.

\item Economic markets: Game theory is being assisted by machine learning (ML) and AI to analyze auctions, pricing strategies, market competition for fairness and efficiency, prediction of market movements and optimization of resource distribution in supply chains. 

\item Cybersecurity: This is certainly a timely area of application. IGT efficiently models interactions between attackers and defenders, with optimal defense strategies  against potential cyber attacks to predict and mitigate threats in intrusion detection systems. Of particular importance is the design of effective protocols and decision-making mechanisms in adversarial security settings, including airport and seaport patrols on the models of algorithmic game-theoretic models for K9-patrols at LAX and national parks protection against poachers.

\item Autonomy: Game-theoretic models are definitely useful for autonomous vehicles to avoid collision, optimize traffic flow, negotiate with other vehicles. Task allocation and path planning problems for multi-robot systems can be solved through intelligent cooperative game theory. Overall, coordination and decision-making of autonomous agents in multi-agent environments greatly benefit from the development of intelligent game theory.
\end{itemize}

\subsection{Testing the Efficiency of Intelligent Game Theory}

\begin{itemize}

\item Computing equilibrium solution and in particular Nash Equilibrium involves some intractable computational complexity. Approximate Nash Equilibrium can be computed using some learning algorithms  such as Q-learning or deep Q-networks with intelligent agents.

\item in Intelligent Cooperative Game Theory, AI is an efficient tool to compute the well-known solution concept of Shapley Value which allocates payout based on contributions, in particular, in distributed systems or collaborative robotics with intelligent agents working together to complete a target task while deciding on fair compensation for each contribution.

\item Computational complexity is inherent to very large scale problems, which necessitates the use of efficient computation methods such as deep learning and neural networks to find optimal strategies for game-theoretic solution in multi-agents environments.

\item Scalability presents another challenge as the number of intelligent agents increases, along with an exponential growth of the complexity of intelligent game-theoretic models.

\item Another critical challenge is of course coordination and cooperation of intelligent with any central authority to achieve specified goals.
\item We should note here also the importance of ethics considerations in all these applications from market design, to autonomous systems to healthcare systems all requiring fairness and transparency.
\end{itemize}

\section{Differential Game Theory}

Classical game theory has been extended to the dynamic systems modeled by differential equations with strategies evolving over time. \textcolor{blue}{Differential Game Theory} (DGT)\footnote{ The introduction of Dynamics Games is due to Rufus Isaacs, 1951, followed with the 1965 book {\it Differential Games: A Mathematical theory with applications to warfare and pursuit, control and optimization,} Other later contributors include A. Merz and Breakwell, T. Basar et al. } deals with the study of dynamic interactions between players in continuous time, and models situations where players have a certain level of control over the rate of change (differential) of certain variables over time, such as the speed of a vehicle or the flow of resources in a network. Therefore such a game theoretical approach is often used in engineering, economics, and ecology to analyze optimal control problems, resource management, and environmental policy. That is, whenever there is a dynamical optimization or optimal control with strategic interaction. For example

\begin{itemize}

\item Aerospace Engineering: The design of optimal guidance and control systems for spacecraft and aircraft will draw from differential game theory taking into account factors such as fuel consumption, aerodynamic performance. and environmental constraints. For example, an aircraft landing subject to wind disturbance. Air traffic control.
    
\item Control Systems and Robotics including autonomous vehicle coordination.

\item Ecology: the study of the dynamics of populations and ecosystems, assuming individuals have control over their reproductive rates or foraging behaviors. The differential game theoretic approach would contribute to identity optimal management strategies for conservation, pest control, or fisheries.
    
\item Military and Defense where there is a need to model competitive strategies between adversaries in conflict situations.

\item Economics and Management Science: Differential Game Theory is effective in analyzing strategic interactions between firms in dynamic market, where production rates, pricing strategies, and investments in research and development are traditionally controlled by the firms. 

\end{itemize}

Differential game theory can be seen indeed as an extension of optimal control where game theory provides equilibrium concepts including Nash equilibrium and Stackelberg equilibrium. See Dockner et al. 2000.
Solution techniques in optimal control include the {\it maximum principle} (Pontryagin et al. 1962) and {\it dynamic programming}. See Bellman's principle of optimality, Belleman et al. 1957. The theory of differential game adds the concept of information structure: {\it open loop} to mean strategies are function of time only, and correspond to the maximum principle, {\it closed loop no memory} shows up in dynamic programming for strategies function of time and the current state, and {\it closed loop memory} for strategies function of time, current and past states. 

The types of differential games are determined by the information patterns $\eta_i(t)$ available to player $P_i$ at time t and the payoff structures: open-loop information or closed loop information.
To illustrate, in a n-players differential game with infinite horizon and state $x$ each player $P_i$ has a control $u_i$ and optimizes an objective function subject to a differential equation. For instance

\begin{equation}
    \max_{y(.)}\int_0^{\infty} e^{-rt} [\beta \gamma (t) -0.5 y^2(t) - 0.5 \gamma s^2(t) ]dt,
\end{equation}
subject to 
$$ \dot{s(t)} = y(t) - \delta s(t),\quad s(0) = s_0 
$$
with $r$ the discount rate and $\delta$ the natural assimilation rate.

\section{Algorithmic Game Theory}

A rapidly developing sub-game theory, {\it Algorithmic Game Theory} (AGT) refers to games involving the use of algorithms or any computational protocol or mechanism for decision-making. That is, the study of game-theoretic problems for the "algorithm" and "computational" perspectives: computing equilibria in an efficient, fast, possibly distributed and centralized manner; thus the question: what is the complexity of computing a Nash Equilibrium? In other words, how long does it take until selfish, however rational, players converge to an equilibrium? Daskalakis, Goldberg, and Papadimitrion, Chen and Deng showed in 2006 that convergence to equilibrium can be prohibitively long. In particular under Nash's fundamental theorem stating that {\it Every finite non-cooperative strategic game of two or more players has a (mixed) Nash equilibrium}.

Algorithmic methods have been used extensively for potential and congestion games to compute their pure Nash equilibria. However, in most experiments the running time has been :exponential". See in reference Nisan et Al. 2007.

The key in AGT is on understanding how computational algorithms interact with game-theoretic frameworks: 

\begin{itemize}

\item Computing Nash equilibria in large scale systems. This can be computational difficult, leading to research on approximation algorithms and efficient mechanisms

\item Driving bidding strategies and optimization in Auction systems\footnote{ A fundamental process of buying and selling goods and services, an auction transfers/exchanges goods and services through taking bids as follows: seller wants to sell an item; bidder $P_i$ has a valuation $v_i$. utility of $P_i$ is $v_i-$price paid, and $0$ if loses auction. submit bid to maximize utility. Seller selects the higher bidder.} (e.g. eBay or Google Ads auctions)

\item Managing traffic routing by selecting the best routes or service allocation strategies

\item Determining the assignments of resources on platforms such as online marketplaces, ride-sharing apps matching drivers with riders in order to reach a computational. equilibrium

\end{itemize}

We now present one of the most prominent applications of Algorithmic Game Theory, {\it Security Games.}

\subsection{Security Games}

Game-theoretic methods are applied to the modeling and analysis of security scenarios, involving mostly attacker and defenders interacting strategically with the goal to optimize security measures while accounting for adversarial behaviors. The scenarios include: cybersecurity, military defense, counter-terrorism, patrolling seaports and airports, protection of critical infrastructures and resources, protecting wildlife and forest form poacher and smugglers, curtailing the illegal flow of weapons, drugs and money across international borders.

The science of security games amounts to designing game-theoretic models to predict and counteract the actions of attackers and optimize the use of available resources for defense. Here a Nash equilibrium is reached whenever neither the attacker nor the defender can improve their payoff by changing their strategy/plan given the strategy/plan put in place by the other, under the rationality assumption.

Mixed strategies are of great importance in this class of games, in particular, for large and complex systems. Actions or plans are chosen or design with certain probabilities, e.g., random patrolling of different routes in order to make it harder for the attacker to predict where security is strongest. Indeed, one of the central questions is how to efficiently allocate finite resources to maximize protection against potential attackers (e.g. security personnel, budget, sensors), in order to balance the cost of security measures with the potential losses from a successful attack. 

We distinguish two main classes:

\begin{itemize}
\item Stackelberg security games\footnote{ First introduce to model leadership and commitment by Kiekintveld et al. } with asymmetric relationship between player and most common in the real-world: Defender commits first to strategy, a course of action before any attack (e.g.,  allocate security resources) and then the attacker reacts to the defender's strategy, by endeavoring to exploit the vulnerabilities in the defensive strategy. 

For example the Stackelberg model has led to the development of new algorithms providing randomized patrolling and inspecting strategies. Examples include ARMOR deployed at LAX in 2007 to randomized checkpoints; IRIS for randomized deployment of US Federal Air Marshals since 2009; PROTECT for randomized patrolling of ports by US Coast Guard (port of Boston since 2011 and port of New York since 2012, spreading to other ports in the nation). PAWS tested by rangers in Uganda for protesting wildlife in national parks. MIDAS for the US Coast Guard for the protection of fisheries.

\item Bayesian security games in which players have totally incomplete information about the other players' strategies or payoffs. Mixed strategies are therefore determinant as defender and attacker consider the probabilities of various types of attack and make decisions based on the belief about the adversary's actions

\end{itemize}

This is to say that the science of security seeks to effectively combine game theory, in particular algorithmic game theory, optimization, strategic thinking oftentimes through "war gaming" and regular drills, to model and solve real-world security issues. This field of research endeavors is crucial in modern days for cybersecurity, military defense, infrastructure protection, counter-terrorism, and strives to leverage the advent of artificial intelligence, machine learning and quantum computing. The most relevant challenges include: complexity, uncertainty, oftentimes due to incomplete information. Dynamics nature of threats requiring constant updates of strategies over time. Ethical implications: for instance, surveillance overreach compromising individuals' privacy, counter-terrorism's impact of civil liberties.

\section{Quantum Game Theory}

Recall that "mixed games" resulted from extending the domain of "pure" payoff profile to include probability distributions. The idea of such extension can be pursued in other directions, including in the non-Archimedean/p-adic and the quantum directions. Game Theory is being extended to include the so-called {\it Quantum games} in which players have access to quantum systems or quantum effects while deciding on their strategies. The most important quantum effects include \textcolor{blue}{superposition} and \textcolor{blue}{entanglement} 
and/or the randomness of quantum measurements. Strategies subjected to these effects cannot be emulated in classical games.  It should be mentioned here some of the pioneering works: first, the 1999 seminal paper by Myers showed that a "quantum player" always wins against a classical player. The Pareto-optimality of a Nash equilibrium compared to the classical Nash equilibrium was showed by Eisert, Wilkens and Lewenstein in 2000 by \textcolor{purple}{quantumizing} a two-person game with both players having access to quantum effects. And according to Dahl and Landsburg, (2011) a quantum payoff resulting from using quantum strategies is at least as great as the payoff using classical strategies

\subsection{Quantum Strategy profile}

Meyer first extended players' pure strategies to include {\it pure quantum strategies}, i.e., the "quantumization" of the game, with the appearance of "new" (almost) optimal Nash Equilibria in terms of quantum strategy profiles. A further extension may be considered to include {\it mixed quantum strategies}, i.e., a probability distribution over the pure quantum strategies. 

\subsubsection{ A Quantum Primer}

{\it Superposition and entanglement} are the most distinguishing features on quantum systems: Superposition refers to a system' ability to appear in multiples states simultaneously whereas entanglement describes a strong (non-classical) correlation between quantum systems with no classical analogue.

We first present some of the key principles of quantum information

\begin{enumerate}

\item Principle 1: A quantum physical system corresponds a \textcolor{blue}{complex Hilbert space $\mathbb{H}$}, that is, a complete Inner Product space $\langle \mathbb{H}, <.,.>\rangle$ over the complex number $\mathbb C$ 

\item Principle 2: The \textcolor{blue}{state $\ket{\psi}$} of a quantum system is given by a unit vector $v\in \mathbb H.$ of the associated Hilbert space.

\item Principle 3: For every \textcolor{blue}{orthonormal basis} $\cal{B}=$$(b_i)_{i=1,\ldots,n}$ in the Hilbert space, there is an associated \textcolor{blue}{ measurement:}\footnote{A measurement can be thought of as an abstract physical process implemented to always return an outcome $\cal{O},$ depending on the state in which the system is. It is probabilistic. } the outcomes of the measurement are in $\cal{B}.$ The probability of outcome $x=\sum_{i=1}^n <b_i,x>b_i \in \cal{B}$ when the system is state $v\in\mathbb H$ is $\|<x,v>\|^2$ ({\it Born rule}). If the outcome of the measurement is $x$ then after measurement the system will be in state 
$\frac{<x,v>}{\|<x,v>\|} x.$ ({\it Wave function collapse})

\item Principle 4: \textcolor{blue}{Time Evolution} $E:\quad S(\mathbb H) \rightarrow S(\mathbb H)$ is linear, with $S(\mathbb H)$ the set of states in $\mathbb H.$ ($E$ is the operator). The linearity implies the preservation of superposition, e.g., $E(a\ket{0}+b\ket{1}) = a E(\ket{0}) + b E(\ket{1}).$ 

\item Principle 5: There is a measurement associated with every \textcolor{blue}{complete family of orthogonal projections.}\footnote{A Projection $P$ is a linear operator with $P^2=P.$ A set of projections $\{P_1,\ldots,P_n\}$ such that $P_iP_j=0$ and $\sum P_i =\mathbf 1$ is called  a complete family of orthogonal projections.}

\item Principle 6: For $\mathbb{H}_1$ and $\mathbb{H}_2$ two Hilbert spaces representing tow quantum physical systems, the Hilbert space of the join system is the \textcolor{blue}{tensor product} 
\textcolor{purple}{$\mathbb{H}_1 \otimes\mathbb{H}_2.$} The above family is a \textcolor{blue}{projective measurement} if the projections are also self-adjoint $P^*=P.$

\end{enumerate}

The quantum state of an isolated system is described by a vector in a  complex Hilbert space $\mathbb{H}(\mathbf(C),$ usually finite dimensional. \textcolor{blue}{qubits} are 2-dimensional systems, and the quantum analog of the classical "bits", represented, using Dirac's famous braket notation
$$
\ket{\psi} = \left[ \begin{array}{c} \alpha\\\beta \end{array} \right],\quad \alpha, \beta \in \mathbf{C}.
$$
The particular states of qubits are
$$
\ket{1}=\left[ \begin{array}{c} 0\\1\end{array} \right]\quad \ket{0}=\left[ \begin{array}{c} 1\\0\end{array} \right]
$$
leading to the general quantum of a qubit as a complex linear combination of these canonical states
$$
\ket{\psi} = \alpha \ket{0} + \beta\ket{1},
$$
that is, a \textcolor{blue}{superposition state}, a defining characteristic of a quantum system with non classical analogue.

We label  $\{\ket{0}, \ket{1} \}$ a set of two normalized, mutually orthogonal states, for qubit, defining these states:
$$
\ket{1}=\left( \begin{array}{c} 0\\1\end{array} \right)\quad \ket{0}=\left( \begin{array}{c} 1\\0\end{array} \right)
$$
The state of the qubit is represented either in a classical or a superposition state
\textcolor{purple}{
\begin{align}
|\ket{\Psi_{classical}} = &\gamma \ket{n},\quad n=0,1\quad |\gamma|^2=1\\
|\ket{\Psi_{superposition}} = & \alpha \ket{0} + \beta \ket{1} \quad |\alpha|^2 + |\beta|^2 =1
\end{align}
}
Note that when observed the superposition collapses in the state in which it is observed.

A system state can be also represented  by positive, semi-definite, Hermitian matrices of trace 1 called \textcolor{blue}{density matrices.}

The probability distribution is {\it separable} when each player can choose their mixed strategy (i.e., through probabilities) independently of each other. However, there are game scenarios in which the probability distribution is not separable, with strategies {\it correlated.} These are the so-called \textcolor{blue}{correlated game.} 

\subsubsection{Quantum Game Play: Quantumizing The Battle of the Sexes}

Below we illustrate using one of well-know classical games, \textcolor{blue}{the Battle of the Sexes.}\footnote{The game first appeared in the book "Games and Decisions" by Duncan Luce and Howard Raiffa in 1957.} 

The game play: A {\it married couple}, Husband (H) and Wife (W) are planning a weekend outing "together" (keyword). Madame prefers to go to the Opera (O) while Monsieur (H) would prefer to go a football game (F). We also assume the following utility increment: 1 if they go to the event of their preference and 0 otherwise. Each gets an increment of 2 for going to  the event together and 0 for going separately. Therefore the game play has: a set of 2 players $\{W, H\}$; a set of 2 {\it pure strategies} $\{F,O \}$. and the following payoff/utility matrix

 \begin{table}[h!]
    \setlength{\extrarowheight}{2pt}
    \begin{tabular}{cc|c|c|}
      & \multicolumn{1}{c}{} & \multicolumn{2}{c}{Husband}\\
      & \multicolumn{1}{c}{} & \multicolumn{1}{c}{$O$}  & \multicolumn{1}{c}{$F$} \\\cline{3-4}
      \multirow{2}*{Wife}  & $O$ & $(3,2)$ & $(0,0)$ \\\cline{3-4}
      & $F$ & $(0,0)$ & $(2,3)$ \\\cline{3-4}
    \end{tabular}
  \end{table}

We then proceed with solving the game. 

\begin{enumerate}
\item This game has clearly two pure Nash Equilibria at $(F,F)$ and $(O,O).$ This prompts the need to randomize.
\item We look for a mixed Nash Equilibrium: The Wife chooses Football  with probability $p$ and the husband with probability $q.,$ leading to the updated mixed matrix


Upon randomization by both Wife and Husband to reach indifference between going to Football and Opera (i.e., equal payoff for going to either event), we then solve to  obtain $p=\frac{2}{5}$ and $q=\frac{3}{5}.$ As these are independent probabilities, the probability that both go to the football game is the product, i.e., $\frac{6}{25}$.

Similarly, we have both at the opera with the same probability $\frac{6}{25},$ but husband at the football while wife at the opera has the probability of $\frac{9}{25}$ and husband at the opera while wife at the football has the probability of $\frac{4}{25}.$

Therefore the mixed Nash Equilibrium $\sigma^*=(\sigma^*_1,\sigma^*_2)$ is given by
$$
(\frac{2}{5}, \frac{3}{5});\quad  (\frac{3}{5}, \frac{2}{5}).
$$

Consequently this game exhibits three Nash Equilibria, consisting of 2 pure and 1 mixed, prompting the concern of how the players should decide. This indicates that simple solutions are not always available even in the simplest game play. The most likely outcome is $(F,O)$ if both decide to randomize with a probability of $\frac{9}{25}.$  Coordination would be a better outcome for the married couple.

\item This game play may be extended to an evolutionary game scenario.

\item We now consider the extension to this game to a quantumization of the play
\end{enumerate}

In the quantum version the players Husband and Wife use quantum strategies that are superposition and entanglement of classical strategies 
The game is quantumized by choosing \textcolor{blue}{quantum qubits systems} given by 2 dimensional Hilbert spaces $\mathbb{H}_w$ and $\mathbb{H}_h$ associated with the players wife and husband. 

Then we consider their {\it tensor product} $\mathbb{H}_w \otimes \mathbb{H}_h$ with it orthonormal basis denoted using the classical pure strategies $O$ and $F$, now considered as vector states $\ket{o}$ and $\ket{F}$

$$
B=(\ket{OO}, \ket{OF}, \ket{OO}, \ket{FO}, \ket{FF}
$$
We proceed using the quantum formalism recalled above in the primer. In particular, we deal with {\it density matrices} associated with the strategies. For example, wife and husband have now at their disposal the following {\it entangled state}
$$
\ket{\psi} = \alpha \ket{OO} + \beta \ket{F},\quad |\alpha|^2 + |\beta|^2 = 1
$$
with the associated density matrix $\rho$
$$
\rho = |\alpha|^2 \ket{OO}\ket{OO} +\alpha\bar{\beta} \ket{OO}\ket{FF} + \bar{\alpha}\beta \ket{FF}\ket{OO} + |\beta|^2\ket{FF}\ket{FF}
$$
In this particular game of the Battle of the Sexes, one exploits the entangled strategies to obtain a \textcolor{blue}{Nash Equilibrium Pareto-optimal} as the unique viable solution of the game (giving an even higher reward to both wife and husband than any other solution): in summary, a quantum strategy exhibit a more attractive solution than the classical strategy. However, it should be noted that only using entangled strategies allow the players wife and husband to reach a unique solution, both getting the same expected payoff. The interpretation is that entanglement of the strategies $O$ and $F$ actually represents a strong correlation, directing wife and husband to act the same way. That is, wife and husband, by entangling their strategies, must now play the same strategy with the option to choose sometimes opera (O) and sometimes football (F), obtaining the best reward playing $O$ one-half of the times and playing one-half $F.$ This is indeed more realistic for a married couple, as if marriage, a strong correlation, is here "quantumized into entanglement"! See more details on the computation in references Marinatto and Weber, 2000, Sowmitra 2023.

\begin{remark}(Commnents)
There is a considerable literature on quantum game. Classical game theory has been applied for various modern decision-making processes from diplomacy to economics to cultural dynamics to national security and to agricultural and biological sciences. For example, the game of chicken and the cuban missile crisis. As a recent extension of the classical game theory. quantum game theory is seeing many practical applications, leveraging the superposition and entanglement features available in quantum games. One can also mention work on quantum routing game.
\end{remark}

\section{Non-Archimedean Game Theory}

{\it Non-Archimedean} refers to the violation of the \textcolor{blue}{Archimedean Principle} \footnote{ Given by Otto Stolz in the 1880s in honor of the ancient Greek Geometer and Physicist Archimedes of Syracuse who stated the principle as Axiom V in his {\it On the Sphere and Cylinder} }. It suffices to verify that there are no \textcolor{blue}{infinitesimals elements} to confirm the Archimedean Principle. In other words, a valued field $<\mathbb K, |.|>$ is \textcolor{blue}{Archimedean} if 

\textcolor{purple}{$$
\forall x\in \mathbf K^*, \exists\quad  n\in\mathbf N \quad |\quad  |\underbrace{x+\ldots+x}_{n-times}|=|nx| > 1.
$$}

Classical game theory is studied in an Archimedean/Euclidean space. Therefore \textcolor{blue}{ Non-Archimedean game theory} is studied using non-archimedean spaces, in particular, through p-adic analysis.\footnote{p-adic numbers, ultrametric spaces, non-archimedean numbers, isosceles spaces all express the same idea. Kurt Hensel initiated the p-adic analysis in 1897.}
Archimedean absolute value $|.|$ satisfies the {\it Triangle Inequality}, whereas the p-adic absolute value  satisfies the \textcolor{blue}{ultrametric (strong) inequality} 
\textcolor{purple}{$$
|x+y|\le \max(|x|,|y|).
$$}

Such an inequality leads to a peculiar geometric space where among other things, all triangles are isosceles, every point in a disc is also a center of the disc, all open balls are also closed (clopen balls).  Here also there are non-constant functions with zero derivatives (called pseudo-constants). The p-adic space is totally disconnected. The non-archimedean or p-adic approach handles beautifully hierarchical and fractal-like structures.
It goes without saying that applying p-adic analysis to game theory requires a strong understanding of p-adic number system and all game theory concepts. 

\subsection{p-adic Primer}

Here is a p-adic refresher sufficient for this compendium.
The set of real numbers $\mathbf R=\mathbb{R}_{\infty}$ is indeed a completion of the rational numbers $\mathbb Q$ with respect to the Euclidean or Archimedean norm denoted here \textcolor{purple}{$|.|_{\infty}$} which satisfies the well-known Triangle Inequality $|x+y|_{\infty}\le |x|_{\infty} + |y|_{\infty}$

From the canonical representation of a positive integer \textcolor{purple}{$ n=\prod_{i=1}^{s}p_i^{k_i}$} and the normalized form of a non-zero rational 
\textcolor{purple}{$ r=p^k\frac{a}{b},$} $p$ prime with $a$ and $b$ co-prime, the \textcolor{blue}{ p-adic absolute} \textcolor{purple}{$|.|_p$} is defined as
\textcolor{purple}{\begin{equation}
|r|_p = \left\{ \begin{array}{rcl} p^{-k} & \mbox{for} & r \ne 0 \\
0  &\mbox{for} & r=0
\end{array} \right.
\end{equation}}

The p-adic absolute value $|.|_p$ or ultra-norm and its induced p-adic distance \textcolor{purple}{$d_p(x,y)=|x-y|_p$} satisfy the {\it Ultrametric or Strong Inequality}

\textcolor{purple}{\begin{align}
&|x+y|_p \le \max (|x|_p, |y|_p)\\
& d_p(x,z) \le \max d_p(x,y) + d_p(y,z)
\end{align}}

The inequality changes to equality for $ |x|_p ne |y|_p$ or $d_p(x,y) \ne d_p(y,z)$
In other words, all triangles in the p-adic space are isosceles.

The completion of $\mathbb Q$ with respect to the p-adic norm $|.|_p$ gives the set of p-adic numbers $\mathbb Q_p$ with its associated series expansion form
\textcolor{purple}{$$
\mathbf Q_p \ni r = \sum_{k=\nu}^{\infty} r_kp^k |\quad r_k\in\{0,1,\ldots p-1\}\subset \mathbf Z.
$$}

The expansion is also symbolically represented using the "decimal" or {\it radix} point (as in the real case) by 
\textcolor{purple}{$r= r_{\nu} \cdots r_{-1}\bullet x_0x_1\cdots x_n\cdots$.}

The set $\{0,1,\ldots p-1\}$ is the set $D_p$ of digits in the representation of p-adic numbers.  The set $\mathbf Q_p$ is the field of fractions of $\mathbf Z_p$ the space of {\it p-adic integers}, the compact unit disk 
$\mathbf Z_p=\{ x\in\mathbf Q_p |\quad |x|_p \le 1 \}$

Other peculiarities on the p-adic analysis include 

\begin{itemize}

\item The set $\mathbb Q_p$ is \textcolor{blue}{locally compact and totally disconnected}: That is, calculus is performed with no expectation of "reasonable" analyticity.
\item While a natural geometric ordering for the usual Euclidean metric is given by the real number line, in the ultrametric case a \textcolor{blue}{hierarchical tree} provides such an ordering.

\item Given a p-adic infinite expansion \textcolor{purple}{$\mathbb{Z}_p \ni x=x_0+\cdots+ x_{n-1}p^{n-1} + O(p^n)$}, we have $ \bar{x} =x_0+\cdots+ x_{n-1}p^{n-1} $ as its p-adic approximation, i.e., $x\equiv \bar{x}\quad mod\quad  p^n.$ 

And $|x-\bar{x}|_p \le p^n$ ensuring convergence. $n$ is the {\it order} or {\it absolute precision} of $\bar{x},$ with the {\it relative} precision given by $n=min\{i\in\mathbf Z, n_i\ne 0\}.$

\item Relevant to game theory, in particular, the computation of Nash Equilibria, is the fact that \textcolor{blue}{{\bf p-adic errors do not add}}: 

That is, 
\textcolor{purple}{$x+O(p^{k_1}) + (y + O(p^{k_1}) = x + y + O(p^{\min(k_1, k_2)})$}. This is a tremendous advantage over the real precision with its compounding of round-off errors. (Gregory, 1980. Krishnamurty, 1977. Dixon, 1982. Ruffa et al. 2016))

\item For $\nu_p(x):= \max\{k\in\mathbf Z: p^k | x$\} called the \textcolor{blue}{p-adic valuation of x,} we also have 
$x+O(p^{k_1}) \times (y + O(p^{k_1}) = x\times y + O(p^{\min(k_1, \nu_p(y)k_2, k_2+\nu_p(x))})$

\item p-adic computation also involves the so-called \textcolor{blue}{Hensel Code} originated from the {Hensel's Lifting Lemma}. (Gouv\^ea, 1997. Katok, 2007, Robert, 2000) 
\end{itemize}

\subsection{Applying p-Adic Analysis to Game Theory}

In light of all the above applying p-adic analysis to game theory is indeed justified. 

\begin{itemize}

\item For one, preferences in decision-making are modeled by utility/payoff functions, which are traditionally real-valued, but can, alternatively, be expressed using p-adic numbers, to account for a different scale of the idea of "closeness" or "distance" between choices or strategies, in particular, in situations where relative versus absolute difference between strategies proves more relevant.

That is, to leverage the ultrameticity of p-adic topology, i.e., using the ultrametric distance between strategies instead of the Triangle Inequality distance. In some instances, p-adic valuations may provide a better working framework for decision-making processes with imperfect information.

\item As recalled, p-adic or ultrametric spaces have peculiar geometry and topology, which could be more suitable to analyze the strategy space in games; indeed, the ways strategies are related to each other certainly influences the equilibrium or any outcome of the game. Note that p-adics are closely tied to local structures, thereby emphasizing local differences in payoffs.

\item Evolutionary Game Theory (EGT) is classically studied in Euclidean spaces. One may now consider a \textcolor{blue}{p-adic evolutionary game theory:} mutations, which play important role in EGT can now be thought as discrete steps appearing according to some p-adic valuations, with impactful changes and dramatic evolutionary shifts resulting from very small mutations in the close vicinity of the strategies.  

\item In EGT, there is a potential of players refining their strategies to improve their relative payoffs. A p-adic evolution could model such refinement process, to include a player's strategy that evolves in a p-adic neighborhood with gradual adjustment to optimal and desirable strategies.

\item Notably, an evolutionarily stable strategy (ESS) can be evaluated through an ultrametric distance instead of the classical Euclidean/Archimedean metric, recalling that p-adic closeness differs for the Euclidean one. For example, 2 is 7-adically closer to 51 than it is to 1 because $d_7(2,51)=i/19$ whereas $d_7(1,2)=1$.

Or 3-adically one can express $-1$ as $-1=2+2.3+2.3^2+\cdots + 2.3^n+\cdots $

Therefore p-analysis may lead to new evolutionary equilibria with a \textcolor{blue}{Nash Equilibrium.}

\item Recall the hierarchical structures of p-adic spaces. This could be indeed used to describe different layers of strategy evolution, contrasting localized evolutionary changes within smaller subgroups in large populations with the broader evolutionary dynamics of the entire population. This will better describe the real-life interplay  between local and global  evolutionary pressures, For examples, in the socio-cultural-economic scenario below, the interplay between local Nash Equilibria and the central Nash Equilibrium: Small changes in local subpopulations of a highly stratified society, such as the American society, may easily propagate across the entire society, possibly directing or redirecting overall strategy shifts in a highly nonlinear manner (not excluding chaotic one).

\item  \textcolor{blue}{p-Adic Probability Distribution}: In classical game theory and its related EGT, mixed strategies are considered using a classical probability distribution of the state vectors. p-adic analysis leads to a development of a theory of \textcolor{blue}{p-adic probability} taking values in the p-adic numbers.: frequencies, relative or otherwise, take place in the set \textcolor{purple}{$[0,1]\subset\mathbb Q$} whose closure is the same set in the Archimedean topology but whose closure in the p-adic topology is the whole set $\mathbb{Q}_p$ of p-adic numbers. Consequently, all possible sequences in $\mathbb{Q}_p$ summing up to 1 are considered legitimate p-adic probability distributions. For example, a sequence such as \textcolor{purple}{$ \{1,-5, -1, 6\}$} which is certainly not a standard probability distribution. That is, in the p-adic probability values greater than 1 or even negative are allowable.

\end{itemize}

\section{p-Adic Quantum Game}

 The p-adic analysis has received many concrete applications in recent years, in particular, applications to  physical theories; due to the existence of the {\it Planck length $l_p\approx 10^{-35}$} in quantum gravity, many non-Riemannian model were being considered, to include the models based on the non-Archimedean field $\mathbb{Q}_p$ of p-adic numbers, which exhibits a \textcolor{blue}{fractal-like, ultrametric and hierarchical structure}, as recalled above. A particular important application has been in computer science and cryptography pertaining to the generation of the so-called {\it pseudorandom numbers} and {\it uniform distribution of sequences.} 

 There is also a great interest in developing a p-adic quantum information theory alongside the quantum information theory. However, we present here some preliminary ideas and work on \textcolor{blue}{{\bf p-Adic Quantum Game Theory}}, a theory we are actively developing drawing from all the above combining the quantum and p-aidc primers, with the goal to leverage the distinguishing features of the p-adic analysis and the quantum approach. 
 
 To recap, quantum game features {\it superposition and entanglement} of players' strategies, leading oftentimes to a "better" and more efficient Nash equilibrium as in the case of the quantum Battle of the Sexes.

 On the other hand, there is a tremendous advantage in considering a p-adic analysis of games, in which players' payoff function has now a co-domain in the p-adic number system $\mathbb{Q}_p.$ A p-adic neighborhood of players' strategy accommodates best some gradual adjustment to optimal strategies in evolutionary game theory, as well as the consideration of different layers of strategy time-evolution. Let also mention that one could successfully leverage the \textcolor{blue}{p-adic probability} and a \textcolor{blue}{p-adic time} in evolutionary game theory. Note that in classical evolution theory, time is Euclidean/Archimedean, i.e. classical time is merely statistical (topologically and metrically like the real line),  whereas a \textcolor{blue}{p-adic time} could better describe the \textcolor{blue}{circularity of time periods} (the non-linearity of time viewed as a series of repeating cycles.).

\subsubsection{Tenets of p-adic quantum game}

The formulation of standard quantum theory relies fundamentally on the field $\mathbb C$ of complex numbers, which can also be regarding as a quadratic extension of the real field $\mathbb R.$ The complex conjugation is the natural involution on the complex numbers field. The fundamental requirement of the p-adic quantum game, is to replace the field $\mathbb C$ with a quadratic extension of $\mathbb{Q}_p.$ We follow here the approach taken by Aniello et al. See Aniello, 2023 for more details. That is, we consider a non-quadratic element in \textcolor{purple}{$\mathbb{Q}_p, $ i.e., $\mu \in \mathbf{Q}_p\quad | \quad \mu \notin (\mathbf{Q}_p)^2$} to finally define

\textcolor{purple}{\begin{equation}
\mathbf{Q}_p(\mu):=\{z= x + y\sqrt{\mu}\quad | \quad x, y \in \mathbf{Q}_p\}
\end{equation}}

$\mathbf{Q}_p(\mu)$ is a field extension of $\mathbf{Q}_p,$ a 2-dimensional $\mathbf{Q}_p$ vector space, endowed with a non-Archimedean valuation

\textcolor{purple}{\begin{equation}
    |z|_{p,\mu}:= \sqrt{|z\bar{z}|_p}
\end{equation}}

 with $z= z= x + y\sqrt{\mu},\quad \bar{z} = x -y\sqrt{\mu},\quad z\bar{z} = x^2 -\mu y^2\in \mathbf{Q}_p $

There are many non-isomorphic quadratic extensions of $\mathbf{Q}_p;$ for example, for $\eta =-1$ and prime $p\equiv 3 (mod 4)$ is a non-quadratic element of $\mathbf{Q}_p.$ (e.g., for p=7, take $\mu\in\{-1, 7,-7\}$. For $p=2$ take any $\mu\in\{2,3,5,7,10,14\}.$

Once a quadratic extension has been determined, we then define a \textcolor{blue}{p-adic Hilbert space} $\langle \mathbb{H}_p, \|.\|, ,<.,.>\rangle$ where, $<.,.>$ is a {\it non-Archimedean inner product} and $\|.\|$ is an {\it ultranorm} such that, in general, the ultranorm does not stem directly from the inner product.

Armed with such a definition of p-adic Hilbert space, for example, a \textcolor{blue}{p-adic qubit} is a quantum system described by a 2-dimensional p-adic Hilbert space.

\begin{remark}
\noindent
\begin{itemize}

\item In addition to the fact that for p-adic Hilbert space, the ultranorm is not directly coming from the inner product, i.e., in general we have $\|z\| \ne \sqrt{|<z,z>|},$ the space may contain non-vectors $z$ with $<z,z>=0,$ i.e., $z$ is an {\it isotropic} vector. For example, $z=b_1 + \sqrt{-1} b_2,$ where $b_1,b_2$ are in an orthonormal basis of $\mathbb H.$

\item  Here \textcolor{blue}{p-adic quantum states} are p-adic statistical operators defined as linear operator $\rho\in \cal{L}(\mathbb H)$ such that $\rho=\rho^*$ and $trace(\rho)=tr(\rho) = 1,$ i.e., a \textcolor{blue}{trace one self-adjoint linear operator} forming the set $\cal{S}(\mathbb H).$

\item The description of \textcolor{blue}{p-adic quantum measurement} (i.e., the {\it observables},) is given by self-adjoint operator valued measures (SOVM), a family $\cal M=(M_i)_{i\in I}$ (I a finite index set)  of self-adjoint operators in $\mathbf H$ such that $\sum_I M_i=\mathbf I.$

\item We also define a linear functional $ \omega \in \cal{L}(\mathbf H)$ for any fixed $\rho\in \cal{S}(\mathbf H$ as 
$$
\omega_{\rho}(\sigma):= tr(\rho \sigma) \in \mathbf{Q}_p(\mu)\quad \sigma\in \cal{L}(\mathbf H).
$$
That is, $\omega_{\rho}$ is a {\it normalized involution preserving linear functional} on $\cal{L}(\mathbf H).$
\item Therefore the sequence $\{\omega_{\rho}(M_i) = tr(\rho M_i) \}_{i\in I}$ serves as p-adic probability distribution.

 \end{itemize}   

\end{remark}

 \subsubsection{p-adically quantumizing a game}
 
We summarize by describing a procedure to \textcolor{blue}{p-adically quantumize} a strategic game: we first leverage the salient features of the p-adic probability in the following sense:

\begin{itemize}
    
\item We implement the states of an n-dimensional p-adic quantum system, associated with an n-pure strategies game, say, the 2-strategy game of the Battle of the Sexes $S=\{O,F\}$, representing the options of an Opera outing and an evening football game. Initially, we solve the game to identify two pure Nash Equilibria (PNE) at $(O,O)$ and $(F,F)$, which is not quite helpful for the wife and husband players decision-making. 

\item Then the players randomize the O and F strategies via a classical probability distribution, which yields a single Mixed Nash Equilibrium (MNE) at $\sigma^*=(\sigma_w,\sigma_h),$ with $\sigma^*_w=(p^*,1-p^*),$ and $\sigma^*_h=(q^*,1-q^*).$ Having now three Nash Equilibria is not ideal for the married couple decision-making, which prompts a need of a better coordination; being in a marriage could stand for some type of correlation between the couple's pure strategies. 

\item A strong, non-classical correlation might be required to generate an improved Nash Equilibrium. The pure strategies are therefore quantumized as described above, leading to the superposition and entanglement of the pure strategies; the entangled strategies produce a much better Nash Equilibrium, one that is Pareto-Optimal.

\item The efficiency of this unique Pareto-Optimal Nash Equilibrium may be greatly improved by leveraging the distinguishing features of p-adic analysis: 

\item Expressing the payoffs in the p-adic numbers system $\mathbb{Q}_p$ to allow the ultranorm to estimate, for example, the distance between the payoffs, in particular following the entanglement of the strategies.

\item The evolution of the strategies may have multiple layers which can be analyzed through the hierarchical structure of the p-adic analysis

 \item This process of p-adically quantumizing this particular game of the Battle of Sexes may therefore result in a single Pareto-Optimal Nash Equilibrium with high efficiency.

 \end{itemize}
 
\section{A socio-cultural-economical game theoretic model}

We present an example of how an advanced game-theoretic approach, such as a \textcolor{blue}{p-adic quantum game-theoretic} one,  can be used in the least expected area: the socio-cultural-economical dynamic of a society: modeling the strategic interactions of American using ten $(10)$ pure strategies representing the ten core values in the American society

\subsection{Game Theoretical approach to America cultural evolution}

Culture\footnote{Some view culture as a vehicle for providing generally accepted solutions to problems. Others define it as part of the knowledge owned by a substantial segment of a group but not necessarily by the general population.} consists of \textcolor{blue}{preference distribution} along with \textcolor{blue}{equilibrium behavior distribution}, but also individuals in any given culture are multidimensional. Immigration has been a major driving and defining factor of American culture. The behavior and choices of immigrants are influenced by many factors, e.g., home country traditions, in interaction with other American, while the preferences and equilibrium behaviors are shaped by {\it socialization and self-persuasion} and vary continuously, impacted by cross-cultural contacts, leading as a result, to some type of cultural hybridization. Discrepancies between the ideal and actual choices increase with time and fuel discontent. One would like to make choices in agreement with personal preferences, but the gains and payoffs are dependent upon the degree of choice coordination. The interactions, while seemingly random, are enhanced by factors such as common language, shared symbols and meanings, communication rules, down to culinary habits, family size.

America, as a social entity, has been experiencing rapidly changing demographics, with varying cultural norms and values. An American culture distinctly emerged as a compromise, a {\it modus vivendi} and {\it modus operandi} out of multiple early immigrant cultures, behaviors and preferences. These cross-cultural interactions generate a single hybrid culture, with its own peculiar language and rules of laws. The basic beliefs, assumptions and values by which most American live could be comprehensively described by some ten core values, we here recall, deeply ingrained in most Americans.

\begin{enumerate}

  \item {\bf Individualism and Privacy}. $(A_1):$ Each person is seen as unique, special and independent. Such belief put a premium on individual initiative, expression, orientation and values privacy. Other cultures emphasize group orientation, conformity as a way to societal harmony. American believe they control their own destiny. Other cultures' focus is on extended family, and its corollary of loyalty and responsibility to family, with age given status and respect.
  \item {\bf Equality/Egalitarianism}. $(A_2):$ Each person expects to be treat equally, with a minimum sense of hierarchy, leading to a directness in relations with others, informal sense of self and space. One by-product is to allow the challenging of authority as opposed to respect for authority and social order at all cost prevalent in other cultures. Gender equity rather than different roles for men and women. In other countries people seem to draw a sense of security and certainty from Class and Authority, Rank and Status, considering it reassuring to know, from birth, who you are and where you fit in the complex social system.
  \item {\bf Materialism}. $(A_3):$ Each person enjoys the right to be "well off" and to "pursuit material happiness". People are often judged by their possessions. An undesirable by-product is the well-documented {\it American Greed}. It is a higher priority to obtain, maintain and protect material objects rather than developing and enjoying interpersonal relationships for instance.
  \item {\bf Science and Technology}. $(A_4):$ This value is seen as the driving force for changes and the primary source of goods. Scientific and technological approaches are highly valued, and lead to a problem-solving focus, with a mental process and learning style that are linear, logical and sequential. 
  \item {\bf Progress and Change}. $(A_5):$This leads to a great national optimism and the so-called "Manifest Destiny": nothing is impossible and greatness is ours, as opposed to just accepting life's difficulties. Other societies value stability, tradition, continuity and rich and ancient heritage, seeing change as disruptive and potentially destructive.
  \item {\bf Work and Leisure}. $(A_6):$ This refers to a strong work ethic, and work for its intrinsic values and as the basis of recognition and power. Idleness is seen as a threat to society while leisure should be a reward for hard work and individual achievement. In other cultures, work is seen as a necessity of life, and rewards are based on seniority and relationships.
  \item {\bf Competition}. $(A_7):$ That is, the "Be always First" mentality, encouraging an aggressive and competitive nature. It has led to {\it Free Enterprise} as an economic system, in the belief that competition fosters rapid progress. Cooperation is instead promoted in other cultures.
  \item {\bf Mobility}. $(A_8):$ This refers to a vertical social, economic and physical mobility, in a society with people on the move to better self.
  \item {\bf Volunteerism}. $(A_9):$ Philanthropy is highly valued as a personal choice not a communal expectation, involving associations or denominations rather than kin-groups as in other cultures.
  \item {\bf Action and Achievement Oriented}. $(A_{10}):$ That is, a practical mind set with emphasis on getting things done, a focus of function and pragmatism and a tendency to be brief and business-like through an explicit and direct communication style, with emphasis on content and meaning found in the choice of words. In other cultures the emphasis is on context, meaning is found around the words, and communication is implicit and indirect. Formal handshakes rather than hugs and bows! {\it Dress for success} rather than as a sign of position, wealth, prestige or religious rules.
  \end{enumerate}
Some of these core values, in particular Individualism, Materialism, and Competition,  are perceived sometimes by individuals from other cultural identities, at least during an initial period, with a negative and derogatory connotation. However, a complete cultural integration and economic success in America is for newcomers to be familiarized as quickly and efficient as possible with these values. We consider the cultural evolution strategies to be associated with these ten core values, and we define a probability distribution, Archimedean and Non-Archimedean, for each strategy $A_i,$ $i=1,\ldots,10.$ We are interested in determining the conditions for social stability and progress, and in designing a mathematical framework to understand the American social dynamics, by analyzing the interactive decision-making of all American in relation these ten core values. See Toni, 2017.

All the above aspects of game theoretic approach are applied to decipher the American socio-cultural evolution: 

\begin{itemize}
    \item Americans as players randomize the 10 pure strategies $\{A_1,\ldots,A_{10}\}$, assigning in their daily and random mutual interactions a probability to each one of the core values, leading to an evolutionary matrix game, as described above. The main Nash Equilibrium, represented by the American Constitution, is indeed a mixed Nash Equilibrium co-existing with many local Nash Equilibria, e.g., resulting from strategic interactions e.g., along party lines, socio-economics, religious, ethnic identities, gender, age, etc.... These mixed sub-Nash equilibria are enclosed by limit cycles, we previously called \textcolor{blue}{Nash Limit Cycles}, i.e., self-sustained oscillations in decision-making by American individuals, before converging to or diverging from the equilibrium points.
    
    \item These  10 American pure strategies/core values may also be {\it quantumized}, realistically, using superposition of these values, as well as entanglement of these values. Indeed, these American values are correlated, strongly correlated as times, possibly with some non-classical correlation (\textcolor{blue}{entangled American core values}). This may result in more efficient Nash Equilibrium, Pareto-Optimal, abstract representation of, e.g., a socially optimum American Constitution with necessary amendments for social stability and progress.
    
    \item \textcolor{blue}{p-adically quantumizing} these strategies/core values will certainly lead to an even more optimum Nash Equilibrium leveraging the nuances and hierarchy of the payoff layers. This is indeed the topic of an upcoming research work.
    
\end{itemize}

\section{Concluding Remarks}
This comprehensive compendium on advances in game theory provides a deep insight into the directions of game theory methodologies permeating all aspects of human knowledge.

Machine Learning (ML) and Artificial Intelligence (AI) techniques and their unparalleled computational power are expanding at a lightning speed: OpenAI and its {\it ChatGPT}, the Chinese {\it DeepSeek} are the dominant ones. Algorithmic Game Theory has been aggressively pursuing the design of efficient algorithms for the computing of the Nash Equilibrium. Combining all these techniques with classical, quantum and non-Archimedean/p-adic Game Theory is opening new frontiers in interactive strategic decision-making and predictive analysis across all fields on human knowledge. We are indeed witnessing the emerging of a very powerful tool integrating learning, adaptation and complex decision-making under uncertainty to model and solve real-world problems and beyond: artificial general intelligence, life extension, uploading of consciousness or mind uploading, people and their environment (e.g., food security and biodiversity), and more importantly post-human mathematics (i.e., artificial mathematical creativity). (Ruelle 2013). 
Game Theory, Artificial Intelligence and Quantum Computing are rapidly converging to generate a profound impact on how we understand, address and solve complex strategic problems. 

{}

\end{document}